\documentclass[11pt]{article}
\usepackage[margin=1in]{geometry}
\usepackage{amsfonts,amsmath,amssymb}
\usepackage{algpseudocode}
\usepackage{comment}
\usepackage{enumitem}
\usepackage{subcaption}

\PassOptionsToPackage{obeyspaces}{url}
\usepackage[colorlinks=true,citecolor=blue,urlcolor=blue,linkcolor=blue,bookmarksopen=true]{hyperref}
\usepackage{amsthm}
\usepackage{tikz,ifthen}
\usetikzlibrary{calc}
\tikzset{%
    vertex/.style={draw,circle,ultra thin,white,fill=black},
    edge/.style={ultra thick},
    spedge/.style={line width=6pt,gray!50}, 
    thedge/.style={thin} 
}
\tikzset{/pgf/foreach/parse=true}
\pgfdeclarelayer{back}
\pgfdeclarelayer{fore}
\pgfsetlayers{back,main,fore}

\usepackage{etoolbox}
\patchcmd{\thebibliography}{\leftmargin\labelwidth}{\leftmargin\labelwidth\addtolength\itemsep{-0.1\baselineskip}}{}{}

\author{Emily Heath\thanks{Department of Mathematics, Iowa State University, Ames, IA, USA. \texttt{eheath@iastate.edu}. Supported in part through NSF RTG Grant DMS-1839918.}
    \and Ryan R.\ Martin\thanks{Department of Mathematics, Iowa State University, Ames, IA, USA. \texttt{rymartin@iastate.edu}. Supported in part through Simons Collaboration Grant \#709641.}
    \and Chris Wells\thanks{Department of Mathematics and Statistics, Auburn University, Auburn, AL, USA. \texttt{coc0014@auburn.edu}. This work was primarily done while this author was affiliated with Iowa State University and was supported in part through NSF RTG Grant DMS-1839918.}}

\title{The maximum number of odd cycles in a planar graph}
\date{}

\usepackage[nameinlink,sort]{cleveref}

\newtheorem{theorem}{Theorem}[section]
\newtheorem{lemma}[theorem]{Lemma}
\crefname{lemma}{lemma}{lemmas}

\newtheorem{conj}[theorem]{Conjecture}
\crefname{conj}{conjecture}{conjectures}

\newtheorem{claim}[theorem]{Claim}
\crefname{claim}{claim}{claims}

\newtheorem{prop}[theorem]{Proposition}
\crefname{prop}{proposition}{propositions}

\newtheorem{obs}[theorem]{Observation}
\crefname{obs}{observation}{observations}

\crefname{equation}{equation}{equations}

\theoremstyle{definition}

\newtheorem{defn}[theorem]{Definition}
\crefname{defn}{definition}{definitions}

\crefname{remark}{remark}{remarks}

\crefname{question}{question}{questions}

\crefname{enumi}{property}{properties}

\crefname{equation}{eq.\!}{eqs.\!}
\Crefname{equation}{Eq.\!}{Eqs.\!}

\makeatletter
\DeclareRobustCommand{\crefnosort}[1]{%
    \begingroup\@cref@sortfalse\cref{#1}\endgroup
}
\DeclareRobustCommand{\Crefnosort}[1]{%
    \begingroup\@cref@sortfalse\Cref{#1}\endgroup
}
\makeatother

\newcommand*{\eqdef}{\stackrel{\mbox{\normalfont\tiny{def}}}{=}}        
\newcommand*{\abs}[1]{\lvert #1\rvert}                

\renewcommand*{\epsilon}{\varepsilon}       
\newcommand*{\R}{\mathbb{R}}                
\newcommand*{\mcal}[1]{\mathcal{#1}}        
\newcommand*{\cyc}[2]{\mathbf{G}\bigl(({#2}),C_{#1}\bigr)}

\newcommand*{\tp}[1]{\textproc{#1}}
\newcommand*{\plan}{\mathcal{P}}

\DeclareMathOperator{\numb}{\mathbf{N}}
\DeclareMathOperator{\cp}{\mathcal{C}}

\newcommand*{\doms}{\mathcal{B}}
\newcommand*{\subs}{\mathcal{S}}
\newcommand*{\tum}{\mathcal{T}}
\newcommand*{\nomor}{\subs_\varnothing}
\newcommand*{\onemor}[1]{\subs_{#1}}
\newcommand*{\tumor}[2]{\subs_{{#1}{#2}}}

\DeclareMathOperator{\supp}{supp}   

\begin{document}
\maketitle
\begin{abstract}
    How many copies of a fixed odd cycle, $C_{2m+1}$, can a planar graph contain?
    We answer this question asymptotically for $m\in\{2,3,4\}$ and prove a bound which is tight up to a factor of $3/2$ for all other values of $m$.
    This extends the prior results of Cox--Martin and Lv et al.\ on the analogous question for even cycles.
    Our bounds result from a reduction to the following maximum likelihood question: which probability mass $\mu$ on the edges of some clique maximizes the probability that $m$ edges sampled independently from $\mu$ form either a cycle or a path?
\end{abstract}

\section{Introduction}
For graphs $G$ and $H$, let $\numb(G,H)$ denote the number of (unlabeled, not necessarily induced) copies of $H$ in $G$.
Furthermore, for a planar graph $H$, define
\[
    \numb_\plan(n,H)\eqdef\max\bigl\{\numb(G,H):G\text{ is an $n$-vertex planar graph}\bigr\}.
\]

The study of $\numb_\plan(n,H)$ was initiated by Hakimi and Schmeichel~\cite{hakimi_cycles} who determined both $\numb_\plan(n,C_3)$ and $\numb_\plan(n,C_4)$ precisely.
Later, Alon and Caro~\cite{alon_biclique} continued this line of inquiry by pinning down the value of $\numb_\plan(n,K_{2,k})$ for all values of $k$.
Wormald~\cite{wormald_3conn} and Eppstein~\cite{eppstein_3conn} independently argued that $\numb_\plan(n,H)=\Theta(n)$ when $H$ is a $3$-connected planar graph.
Huynh, Joret and Wood~\cite{huynh_surface}, demonstrated that $\numb_\plan(n,H)=\Theta(n^{f(H)})$ for every planar graph $H$, where $f(H)$ is a graph invariant called the flap number.
See also~\cite{liu} for a further generalization of this result.

Since the order of magnitude of $\numb_\plan(n,H)$ is now understood, the next question is to pin down the coefficient in front of the leading term.
This leading coefficient has been found for several small graphs beyond those mentioned above: $C_5$~\cite{gyori_c5}, $P_4$~\cite{gyori_p4}, $P_5$~\cite{ghosh_planarp5} and $P_7$~\cite{cox_paths}.

The strongest result along these lines to date is that $\numb_\plan(n,C_{2m})=\bigl({n\over m}\bigr)^m+o(n^m)$ for all $m\geq 3$, which was proved for small $m$ in~\cite{cox_paths,cox_cycles} and then extended to all $m$ in~\cite{lv_cycles}.
This paper is motivated by a desire to understand the maximum number of copies of an odd cycle a planar graph can hold.
\medskip

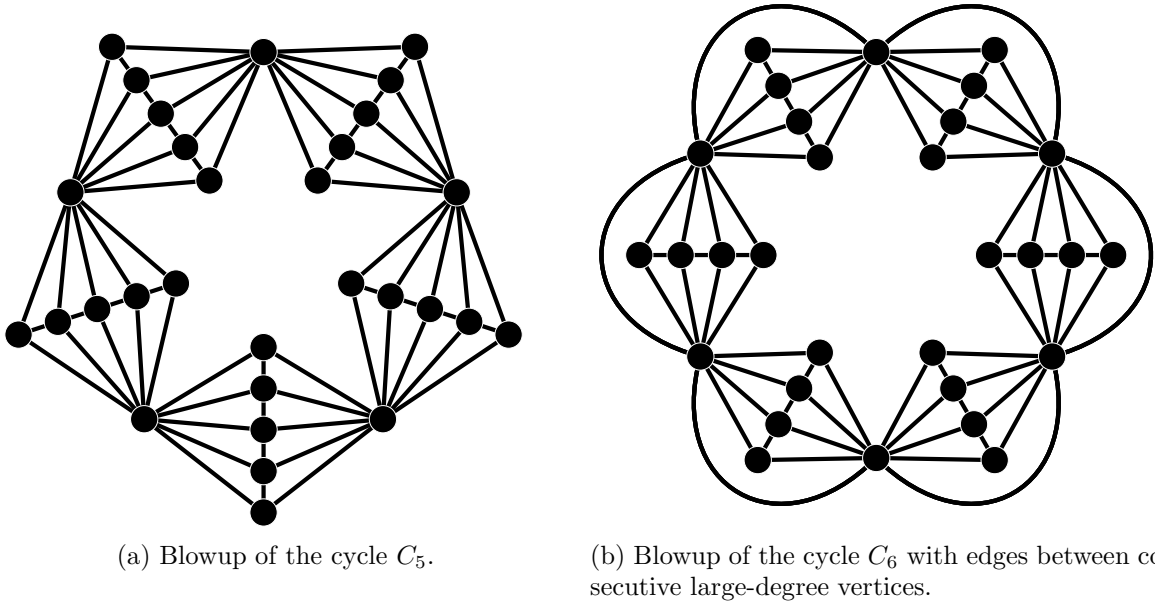
\begin{figure}[ht]
    \centering
    \def\rad{2.7}
    \def\inter{0.55}
    \def\cornera{(-3.8,-3.6)}
    \def\cornerb{(3.8,3.5)}
    \begin{subfigure}[b]{0.48\textwidth}
        \begin{tikzpicture}
            \useasboundingbox \cornera rectangle \cornerb;
            \def\m{5};
            \def\b{5};
            \pgfmathsetmacro{\innrad}{\rad-((\b)/2)*\inter-0.1};
            \pgfmathsetmacro{\outrad}{\rad+((\b-2)/2)*\inter-0.1};
            \pgfmathsetmacro{\tilt}{18};
            \foreach\i in {0,...,{\m-1}}
                {
                    \coordinate (doms\i) at ({\i*360/\m+\tilt}:\rad) {};
                    \foreach\j in {0,...,{\b-1}}
                    {
                        \coordinate (subs\i\j) at ({\i*360/\m+180/\m+\tilt}:{\innrad+\inter*\j});
                    }
                }
            \begin{pgfonlayer}{fore}
                \foreach\i in {0,...,{\m-1}}
                    \node[vertex] at (doms\i) {};
            \end{pgfonlayer}
            \begin{pgfonlayer}{main}
                \foreach\i in {0,...,{\m-1}}
                    \foreach\j in {0,...,{\b-1}}
                        \node[vertex] at (subs\i\j) {};
            \end{pgfonlayer}
            \begin{pgfonlayer}{back}
                \foreach\i in {0,...,{\m-1}}
                    {
                        \pgfmathsetmacro{\next}{int(mod(\i+1,\m))};
                        \foreach\j in {0,...,{\b-2}}
                            {
                                \pgfmathsetmacro{\nextt}{\j+1};
                               \draw[edge] (doms\i) -- (subs\i\j) -- (doms\next);
                                \draw[edge] (subs\i\j) -- (subs\i\nextt);
                            }
                            \pgfmathsetmacro{\nextt}{\b-1};
                            \draw[edge] (doms\i) -- (subs\i\nextt) -- (doms\next);
                    }
            \end{pgfonlayer}
        \end{tikzpicture}
        \caption{Blowup of the cycle $C_5$.~\\~\label{C5blowup}}
    \end{subfigure}\hfil%
    \begin{subfigure}[b]{0.48\textwidth}
        \begin{tikzpicture}
            \useasboundingbox \cornera rectangle \cornerb;
            \def\m{6};
            \def\b{4};
            \pgfmathsetmacro{\innrad}{\rad-((\b)/2)*\inter-0.1};
            \pgfmathsetmacro{\outrad}{\rad+((\b-2)/2)*\inter-0.1};
            \pgfmathsetmacro{\tilt}{30};
            \foreach\i in {0,...,{\m-1}}
                {
                    \coordinate (doms\i) at ({\i*360/\m+\tilt}:\rad) {};
                    \foreach\j in {0,...,{\b-1}}
                    {
                        \coordinate (subs\i\j) at ({\i*360/\m+180/\m+\tilt}:{\innrad+\inter*\j});
                    }
                }
            \begin{pgfonlayer}{fore}
                \foreach\i in {0,...,{\m-1}}
                    \node[vertex] at (doms\i) {};
            \end{pgfonlayer}
            \begin{pgfonlayer}{main}
                \foreach\i in {0,...,{\m-1}}
                    \foreach\j in {0,...,{\b-1}}
                        \node[vertex] at (subs\i\j) {};
            \end{pgfonlayer}
            \begin{pgfonlayer}{back}
                \foreach\i in {0,...,{\m-1}}
                    {
                        \pgfmathsetmacro{\next}{int(mod(\i+1,\m))};
                        \foreach\j in {0,...,{\b-2}}
                            {
                                \pgfmathsetmacro{\nextt}{\j+1};
                               \draw[edge] (doms\i) -- (subs\i\j) -- (doms\next);
                                \draw[edge] (subs\i\j) -- (subs\i\nextt);
                                \pgfmathsetmacro{\distnc}{1.1*(\outrad-\innrad)};
                                \draw[edge] (doms\i) to[out={\i*360/\m+\tilt+45}, in={\i*360/\m+360/\m+\tilt-45}, distance=\distnc cm] (doms\next);
                            }
                            \pgfmathsetmacro{\nextt}{\b-1};
                            \draw[edge] (doms\i) -- (subs\i\nextt) -- (doms\next);
                    }
            \end{pgfonlayer}
        \end{tikzpicture}
        \caption{Blowup of the cycle $C_6$ with edges between consecutive large-degree vertices.\label{C6blowup}}
    \end{subfigure}
    \caption{Two configurations with a large number of copies of $C_{11}$. \Cref{C6blowup} has approximately $n^5/6^4$ such copies, but \Cref{C5blowup} has approximately $2n^5/5^4$ such copies.}
\end{figure}

For $m\geq 3$, a lower bound of $\numb_\plan(n,C_{2m+1})\geq 2m\bigl({n\over m}\bigr)^m-O(n^{m-1})$ is realized by starting with a copy of $C_m$ and replacing each edge $xy$ by a path on approximately $n/m-1$ many vertices and connecting each of these new vertices to both $x$ and $y$ (see~\Cref{C5blowup}).
We believe this construction to be asymptotically tight for all $m\geq 3$, and we make strides toward proving this to be the case.

\begin{theorem}\label{oddcycles}
    The following asymptotic bounds hold:
    \begin{align*}
        \numb_\plan(n,C_5) &= 2n^2+O(n^{2-1/5}), && \text{and}\\
        \numb_\plan(n,C_{2m+1}) &= 2m\biggl({n\over m}\biggr)^m + O(n^{m-1/5}), &&\text{for }m\in\{3,4\}\text{, and}\\
        \numb_\plan(n,C_{2m+1}) &\leq 3m\biggl({n\over m}\biggr)^m + O(n^{m-1/5}), &&\text{for }m\geq 5.
    \end{align*}
\end{theorem}
We note that the constant $3$ can be improved without much effort, especially for larger values of $m$; however, bringing this constant all the way down to our conjectured value of $2$ is currently beyond our reach.

We additionally note that the value of $\numb_\plan(n,C_5)$ was previously determined exactly for all $n$ by Gy\H{o}ri et al.~\cite{gyori_c5} through wildly different means.
However, we include this (weaker) result to demonstrate the method developed in this paper.

\subsection{Reduction to maximum likelihood estimator problems on graphs}

Very generally, maximum likelihood estimator question asks: which probability distribution maximizes the probability of a certain set of observations?
Historically, these questions were focused on determining the member of a family of probability distributions (e.g.\ the family of normal distributions) that best fits a set of observed data.
Recently, Cox and Martin~\cite{cox_paths} showed that bounding $\numb_\plan(n,H)$, assuming $H$ has a special subdivision structure, can be reduced to a question asking: which probability distribution $\mu$ on the edges of a clique maximizes the probability that $e(H')$ many edges sampled independently from $\mu$ yields a copy of $H'$?
While significantly different in scope, this question can be viewed as a ``maximum likelihood estimator question on graphs'' and appears to be absent from the literature, save the papers resulting from this line of inquiry~\cite{antonir_paths,cox_paths,cox_cycles,lv_cycles}.

The biggest success of this reduction to maximum likelihood estimators to date is the proof that $\numb_\plan(n,C_{2m})=\bigl({n\over m}\bigr)^m+o(n^m)$ for every $m\geq 3$, in which case the corresponding maximum likelihood question was solved for small $m$ in~\cite{cox_paths,cox_cycles} and for all $m$ by Lv et al.~\cite{lv_cycles}.
We discuss the actual maximum likelihood question in this case shortly (see \Cref{evenreduction}).

The key contribution of this manuscript is an extension of the methods of Cox and Martin in order to relate the problem of bounding $\numb_\plan(n,C_{2m+1})$ to a maximum likelihood estimator question on graphs.
\medskip

\begin{defn}
    An \emph{edge probability measure} $\mu$ is a probability measure on the edges of some complete graph.
    For a complete graph $K$, we denote by $\Delta^K$ the set of all edge probability measures on $K$.
\end{defn}
For any clique $K$, note that $\Delta^K$ is naturally identified with the (${\abs{V(K)}\choose 2}-1$)-dimensional simplex.
\medskip

For graphs $G,H$, we denote by $\cp(G,H)$ the set of all (unlabeled, not necessarily induced) copies of $H$ within $G$.
Therefore, $\numb(G,H)=\abs{\cp(G,H)}$.

\begin{defn}
    Fix an edge probability measure $\mu\in\Delta^K$ for some clique $K$.
    \begin{itemize}
        \item For a subgraph $H\subseteq K$, define
            \[
                \mu(H)\eqdef\prod_{e\in E(H)}\mu(e).
            \]
        \item For a graph $H$, define
            \[
                \beta(\mu;H)\eqdef\sum_{H'\in\cp(K,H)}\mu(H').
            \]
    \end{itemize}
\end{defn}

$\beta(\mu,H)$ can be viewed as the probability that $e(H)$ many edges sampled independently from $\mu$ form a copy of $H$.

\medskip

The following is one of the key reduction lemmas of Cox and Martin.
\begin{lemma}[Reduction lemma for even cycles~{\cite[Lemma 2.5]{cox_paths}}]\label{evenreduction}
    For every $n$-vertex planar graph $G$, there is an edge probability measure $\mu$ such that
    \begin{align*}
        \numb(G,C_4) &\leq \biggl({1\over 2}\sum_{e\in\supp\mu}\mu(e)^2\biggr) n^2+O(n^{2-1/5}), && \text{and}\\
        \numb(G,C_{2m}) &\leq \beta(\mu;C_m)\cdot n^m+O(n^{m-1/5}), && \text{for all $m\geq 3$},
    \end{align*}
    where the implicit constant in the big-oh notation depends on $m$.
\end{lemma}

The key contribution of this paper is an analogous reduction lemma for odd cycles.
\begin{lemma}[Reduction lemma for odd cycles]\label{reduction}
    For every $n$-vertex planar graph $G$, there is an edge probability measure $\mu$ such that
    \begin{align*}
        \numb(G,C_5) &\leq \biggl(2\sum_{e\in\supp\mu}\mu(e)^2+\beta(\mu;P_3)\biggr)n^2+O(n^{2-1/5}), && \text{and} \\
        \numb(G,C_{2m+1}) &\leq  \bigl(2m\cdot \beta(\mu;C_m)+\beta(\mu;P_{m+1})\bigr)n^m+O(n^{m-1/5}), && \text{for all $m\geq 3$},
    \end{align*}
    where the implicit constant in the big-oh notation depends on $m$.
\end{lemma}

The reduction lemma for even cycles is more general than stated; it actually applies to the class of graphs with linearly many edges and no copy of $K_{3,t}$ for some $t$.
This includes, in particular, the class of graphs embeddable onto surfaces of any fixed genus.
However, the reduction lemma for odd cycles developed in this paper relies critically on the topology of the plane.
\medskip

The majority of this manuscript is dedicated to proving the reduction lemma for odd cycles.
After proving the reduction lemma, we then bound the resulting maximum likelihood questions in order to produce \Cref{oddcycles}.

\begin{theorem}\label{maxlikelihood}
    \begin{align*}
        \sup_{\mu}\biggl(2\sum_{e\in\supp\mu}\mu(e)^2+\beta(\mu;P_3)\biggr) &= 2, && \text{and}\\
        \sup_{\mu}\bigl(2m\cdot \beta(\mu;C_m)+\beta(\mu;P_{m+1})\bigr) &= \frac{2}{m^{m-1}}, && \mbox{for $m\in\{3,4\}$, and} \\
        \sup_{\mu}\bigl(2m\cdot \beta(\mu;C_m)+\beta(\mu;P_{m+1})\bigr) &<  \frac{2.7}{m^{m-1}}, &&\mbox{for all $m\geq 5$.}
    \end{align*}
\end{theorem}

As mentioned previously, the constant $2.7$ can be lowered, especially for larger values of $m$, but it is currently beyond our reach to bring it all the way down to our conjectured value of $2$.
Furthermore, if one seeks only a bound of the form $C/m^{m-1}$ for $m\geq 5$ where $C$ is some absolute constant, then one can na\"ively use the known bounds of $\beta(\mu;C_m)\leq {1\over m^m}$ (Lv et al.~\cite{lv_cycles}) and $\beta(\mu;P_{m+1})\leq {20\over (m+1)^{m-1}}$ (Antonir and Shapira~\cite{antonir_paths}).

The proof of \Cref{maxlikelihood} can be found in \Cref{maxlikelihood-proof}.
There, the three stated bounds are proved separately as \Cref{C5-max}, \Cref{exact34} and \Cref{somebound}, respectively.

\subsection{Notation}
In this paper, all graphs are simple and we use standard graph theory definitions and notation (generally following \cite{west_graph}).
For a graph $G$, we write $e(G)\eqdef\abs{E(G)}$ and $v(G)\eqdef\abs{V(G)}$.
For a vertex $v\in V(G)$, $N(v)\eqdef\{u : uv\in E(G)\}$ denotes the neighborhood of $v$ in $G$.
For disjoint subsets $A,B\subseteq V(G)$, $G[A]$ denotes the subgraph of $G$ induced on $A$ and $G[A,B]$ denotes bipartite subgraph of $G$ induced between $A,B$.


For distinct elements $x,y$, we abbreviate the set $\{x,y\}$ to $xy$, mirroring common shorthand for an edge in a graph.

Throughout this paper, we fix the value $m$ and obtain results to compute upper bounds for $\numb(G,C_{2m+1})$ for graphs $G$ of large order.
As such, all implicit constants in any big-oh notation will depend on $m$.

\section{Preliminaries}



The following fact is well-known:
\begin{prop}\label{denseplanar}
    If $G$ is a planar graph, then $e(G)\leq 3v(G)$,
    If $G$ is a planar bipartite graph, then $e(G)\leq 2v(G)$.
\end{prop}


One immediate consequence of these bounds is that a planar graph cannot have too many vertices of large degree:
\begin{prop}\label{deg-size}
    An $n$-vertex planar graph contains at most $6n/d$ many vertices of degree $\geq d$.
\end{prop}

A considerable number of arguments in this manuscript rely on the fact that planar bipartite graphs are sparse.
In particular, a planar bipartite graph cannot have many vertices of degree $3$ or larger:
\begin{prop}\label{bipartite}
    Let $G$ be a planar graph and fix any $B\subseteq V(G)$.
    If $A\subseteq V(G)\setminus B$ is the set of all vertices $v$ with $\abs{N(v)\cap B}\geq k\geq 3$, then
    \[
        \abs A\leq {2\abs B\over k-2}.
    \]
\end{prop}
\begin{proof}
    Consider the bipartite subgraph of $G$ with parts $A$ and $B$.
    By construction, this subgraph has at least $k\abs A$ many edges.
    Additionally, this is a bipartite planar graph and so it has at most $2(\abs A + \abs B)$ many edges (\Cref{denseplanar}).
    Therefore,
    \[
        k\abs A\leq 2(\abs A+\abs B)\quad\implies\quad \abs A\leq {2\abs B\over k-2}.\qedhere
    \]
\end{proof}

One particular consequence of the above proposition is that a bounded-degree planar bipartite graph cannot contain too many copies of $P_3$:
\begin{prop}\label{forks}
    Let $G$ be a planar bipartite graph with parts $A,B$.
    The number of copies of $P_3$ with both endpoints in $B$ and midpoint in $A$ is bounded above by $\abs A+4d\abs B$ where $d=\max_{v\in A}\deg v$.
\end{prop}
\begin{proof}
    For each positive integer $k$, let $A_k\subseteq A$ denote the set of vertices of $A$ with exactly $k$ neighbors in $B$.
    Then the number of copies of $P_3$ of the desired type is precisely
    \begin{align*}
        \sum_{k=2}^d{k\choose 2}\abs{A_k} &= \sum_{k=2}^d\sum_{i=1}^{k-1}i\abs{A_k}=\sum_{i=1}^{d-1}i\sum_{k=i+1}^{d}\abs{A_k} = \abs{A_2}+\biggl\lvert\bigcup_{k\geq 3}A_k\biggr\rvert+\sum_{i=2}^{d-1} i\cdot\biggl\lvert\bigcup_{k\geq i+1}A_k\biggr\rvert.
    \end{align*}
    By then applying \Cref{bipartite}, we continue to bound
    \begin{align*}
        \abs{A_2}+2\abs B+\sum_{i=2}^{d-1}i\cdot{2\abs{B}\over i-1} &\leq \abs{A_2}+2\abs B+\sum_{i=2}^{d-1}4\abs B \leq \abs{A_2}+4d\abs{B}\leq \abs A+4d\abs B.\qedhere
    \end{align*}
\end{proof}

Finally, we will rely on known orders of magnitude of the maximum number of paths and cycles in a planar graph:
\begin{theorem}[Gy\H{o}ri et al.~\cite{gyori_turan}]\label{pathasymptotics}
    \begin{align*}
        \numb_\plan(n,P_{2m}) &=\Theta(n^m), && \text{for all $m\geq 0$, and}\\
        \numb_\plan(n,P_{2m+1}) &=\Theta(n^{m+1}), && \text{for all $m\geq 0$, and}\\
        \numb_\plan(n,C_{2m+1}) &=\Theta(n^m), && \text{for all }m\geq 1.
    \end{align*}
\end{theorem}
Note that these formulas work even in the trivial cases of $P_0,P_1,P_2$, where $P_0$ is the null graph.
We will need this result, even in the trivial cases, in our proofs.




\section{Proof of \Cref{reduction}: The reduction lemma for odd cycles}
\label{reduction-proof}

Given a planar graph $G$ on $n$ vertices, we will find graphs $G_1$ and then $G_2$ so that the total number of copies of $C_{2m+1}$ in $G$ is the same as in $G_2$, up to a small error term, where $G_2$ will be highly structured and in which counting the cycles is asymptotically equivalent to solving a maximum likelihood problem.

Both $G_1$ and $G_2$ will be so-called tumor graphs:
\begin{defn}[Tumor graph]
    A \emph{tumor graph} is a triple $(H;\doms,\subs)$ where
    \begin{itemize}
        \item $H$ is a graph, and
        \item $V(H)=\doms\sqcup\subs$, and
        \item Every vertex in $\subs$ has at most two neighbors in $\doms$.
    \end{itemize}
    For a tumor graph $(H;\doms,\subs)$, define
    \begin{align*}
        \nomor(H) &\eqdef \{s\in\subs:N(s)\cap\doms=\varnothing\}, &&\text{and} \\
        \onemor x (H)&\eqdef \{s\in\subs:N(s)\cap\doms=\{x\}\},     && \text{for all $x\in\doms$, and}\\
        \tumor xy (H)&\eqdef \{s\in\subs: N(s)\cap\doms=\{x,y\}\},  && \text{for all } \{x,y\}\in\binom{\doms}{2}.
    \end{align*}
    The sets $\subs_{xy}(H)$ are called \emph{tumors}.
\end{defn}
When the tumor graph is understood, we drop the parenthetical and simply write $\nomor,\onemor x,\tumor xy$.

By the definition of a tumor graph,
\[
    \subs = \nomor\sqcup\bigsqcup_{x\in\doms}\onemor x \sqcup\bigsqcup_{xy\in{\doms\choose 2}}\tumor xy.
\]

Notice that our conjectured asymptotic extremal examples for $\numb_\plan(n,C_{2m+1})$ are all tumor graphs (see~\Cref{C5blowup}).
\medskip

We first extract a tumor graph $G_1$ from the original graph $G$.
\begin{lemma}\label{tumorsamecycles}
    Let $G$ be a planar graph on $n$ vertices.
    There is a spanning subgraph $G_1$ of $G$ and a partition $V(G)=\doms\sqcup\subs$ with the following properties:
    \begin{itemize}
        \item $(G_1;\doms,\subs)$ is a tumor graph, and
        \item $\abs\doms \leq O(n^{4/5})$, and
        \item $\deg v\leq n^{1/5}$ for each $v\in\subs$, and
        \item $\numb(G,C_{2m+1})\leq \numb(G_1,C_{2m+1})+O(n^{m-1/5})$.
    \end{itemize}
\end{lemma}

The proof of \Cref{tumorsamecycles} is in \Cref{sec:tumorsamecycles}.
\medskip

The next step will be to prove that, in $G_1$, almost all of the copies of $C_{2m+1}$ alternate between vertices in $\doms$ and vertices in $\subs$ as much as possible.

\begin{defn}[Good cycle]
    Let $(H;\doms,\subs)$ be a tumor graph and let $m$ be a positive integer.
    A copy of $C_{2m+1}$ in $H$ is said to be \emph{good} if it has the form
    \[
        \underbrace{\doms\subs\cdots\doms\subs}_{2m}\subs \qquad\text{or}\qquad\underbrace{\doms\subs\cdots\doms\subs}_{2m}\doms
    \]
    That is, a good cycle contains at least $m$ vertices from $\doms$ and at least $m$ from $\subs$ such that all but one consecutive pair alternates between $\doms$ and $\subs$.

    The number of good copies of $C_{2m+1}$ in $(H;\doms,\subs)$ is denoted by $\cyc{2m+1}{H;\doms,\subs}$.
\end{defn}

Observe that in \Cref{C5blowup}, every good copy of $C_{11}$ has the form $\doms\subs\cdots\doms\subs\subs$, whereas in \Cref{C6blowup}, every good copy of $C_{11}$ has the form $\doms\subs\cdots\doms\subs\doms$.

\begin{lemma}\label{fewbadcycles}
    Let $(G_1;\doms,\subs)$ be a planar tumor graph on $n$ vertices.
    If $d$ denotes the largest degree of a vertex in $\subs$, then
    \[
        \numb(G_1,C_{2m+1})\leq\cyc{2m+1}{G_1;\doms,\subs}+O\bigl(d^3\abs\subs n^{m-2}+\abs\doms n^{m-1}+\abs\subs^3 n^{m-4}\bigr).
    \]
\end{lemma}

By our choice of $G_1$ from \Cref{tumorsamecycles}, $d\leq n^{1/5}$, $\abs \subs \leq n$ and $\abs \doms \leq O(n^{4/5})$.
Thus, $\numb(G_1,C_{2m+1})\leq\cyc{2m+1}{G_1;\doms,\subs}+O(n^{m-1/5})$. The proof of \Cref{fewbadcycles} is in \Cref{sec:fewbadcycles}.
\medskip

Finally, we will find a planar tumor graph $G_2$ which is more refined than $G_1$ but which contains asymptotically the same number of good copies of $C_{2m+1}$ as $G_1$.

\begin{defn}[Benign tumor graph]
    A tumor graph $(H;\doms,\subs)$ is said to be \emph{benign} if whenever $uv\in E(H[\subs])$, then $u,v\in\tumor xy$ for some $xy\in{\doms\choose 2}$.
\end{defn}

Note that our conjectured asymptotic extremal examples for $\numb_\plan(n,C_{2m+1})$ are all benign tumor graphs (see~\Cref{C5blowup}).

\begin{lemma}[Cleaning lemma]\label{cleaning}
    Let $(G_1;\doms,\subs)$ be a planar tumor graph on $n$ vertices.
    There is another planar tumor graph $(G_2;\doms',\subs')$ such that
    \begin{enumerate}[label=(\roman*)]
        \item $(G_2;\doms',\subs')$ is benign, and
        \item $\abs{\doms'}\leq 3\abs\doms$, and
        \item $\cyc{2m+1}{G_1;\doms',\subs'}\leq\cyc{2m+1}{G_2;\doms,\subs}+O(\abs\doms n^{m-1})$.
    \end{enumerate}
\end{lemma}

Because $\abs \doms \leq O(n^{4/5})$ in our case, we would have $\cyc{2m+1}{G_1;\doms,\subs}\leq\cyc{2m+1}{G_2;\doms',\subs'}+O(n^{m-1/5})$.
The proof of \Cref{cleaning} is in \Cref{sec:cleaning}.
\medskip

Putting together the definitions of $G_1$ and $G_2$ along with \Cref{tumorsamecycles,fewbadcycles,cleaning}, we obtain
\begin{align}
    \numb(G,C_{2m+1}) &\leq \numb(G_1,C_{2m+1})+O(n^{m-1/5}) \nonumber\\
    &\leq \cyc{2m+1}{G_1;\doms,\subs}+O(n^{m-1/5}) \nonumber\\
    &\leq \cyc{2m+1}{G_2;\doms',\subs'}+O(n^{m-1/5}) \label{eqn:reduce}
\end{align}

The final step in the proof of the reduction lemma is to actually count good cycles in a benign planar tumor graph.
\begin{lemma}\label{reducetomax}
    If $(G_2;\doms,\subs)$ is a benign planar tumor graph on $n$ vertices, then there exists an edge probability measure $\mu$ on the clique with vertex set $\doms$ such that
    \begin{align*}
        \cyc{5}{G_2;\doms,\subs} &\leq \biggl(2\sum_{e\in\supp\mu}\mu(e)^2+\beta(\mu;P_3)\biggr)\cdot n^2 && \text{and}\\
        \cyc{2m+1}{G_2;\doms,\subs} &\leq \bigl(2m\cdot \beta(\mu;C_m)+\beta(\mu;P_{m+1})\bigr)\cdot n^m &&\text{for }m\geq 3.
    \end{align*}
\end{lemma}
The proof of \Cref{reducetomax} is in \Cref{sec:reducetomax}.
\medskip

In light of \cref{eqn:reduce}, this will complete the proof of \Cref{reduction}.

\subsection{Proof of \Cref{tumorsamecycles}: Finding a tumor graph within a planar graph}
\label{sec:tumorsamecycles}

    Set $d=n^{1/5}$ and define the partition
    \[
        \subs\eqdef\{v:\deg v< d\},\qquad\text{and}\qquad\doms\eqdef\{v:\deg v\geq d\}.
    \]
    By definition, $\deg v\leq n^{1/5}$ for each $v\in\subs$ and $\abs\doms\leq 6n/d=O(n^{4/5})$ by \Cref{deg-size}.

    Now, set $D=n^{2/5}$ and define
    \[
        \doms_{< D}\eqdef\{v\in \doms:\deg v< D\},\qquad\text{and}\qquad \doms_{\geq D}\eqdef\{v\in \doms:\deg v\geq D\}.
    \]
    To begin, let $\subs'\subseteq \subs$ denote the set of all vertices in $\doms$ which have at least three neighbors in $\doms_{\geq D}$.
    Of course, $\abs{\doms_{\geq D}}\leq 6n/D$ thanks to \Cref{deg-size} and so \Cref{bipartite} implies that
    \[
        \abs{\subs'}\leq 2\abs{\doms_{\geq D}}\leq 12n/D
    \]
    as well.
    By \Cref{denseplanar},
    \begin{equation}\label{eqn:smallbigedges}
        e\bigl(G[\subs',\doms_{\geq D}]\bigr)\leq 2\abs{\subs'}+2\abs{\doms_{\geq D}}\leq O(n/D).
    \end{equation}

    Let $G'$ be the subgraph of $G$ formed by deleting all edges between $\subs'$ and $\doms_{\geq D}$.
    After these edges are deleted, every vertex in $\subs$ is now adjacent, in the resulting $G'$, to at most two vertices in $\doms_{\geq D}$.

    \begin{claim}
        With $G'$ defined as above,  $\numb(G,C_{2m+1})\leq\numb(G',C_{2m+1})+O(n^{m-1/5})$.
    \end{claim}
    \begin{proof}
        We simply need to bound the number of copies of $C_{2m+1}$ in $G$ which use some edge within $E\bigl(G[\subs',\doms_{\geq D}]\bigr)$.
        In other words, we need to bound the number of copies of $C_{2m+1}$ of the form $(v_1,\dots,v_{2m+1})$ where $v_1\in \doms_{\geq D}$ and $v_2\in \subs'$.
        \Cref{eqn:smallbigedges} implies that there are at most $O(n/D)$  choices for the pair $(v_1,v_2)$; then there are at most $d$  choices for $v_3$ since $v_2\in S'\subseteq \subs$.
        In total, there are at most $O(dn/D)=O(n^{4/5})$  choices for the triple $(v_1,v_2,v_3)$.
        Finally, by \Cref{pathasymptotics}, there are at most $2\numb(G,P_{2(m-1)})\leq O(n^{m-1})$  choices for the path $(v_4,\dots,v_{2m+1})$.
        This yields a total of at most $O(n^{4/5}\cdot n^{m-1})=O(n^{m-1/5})$  copies of $C_{2m+1}$ which exist in $G$ but not in $G'$, thus finishing the claim.
    \end{proof}
    We may now disregard the graph $G$ and work solely with the graph $G'$.

    Let $\subs''\subseteq \subs$ be the set of all vertices in $\subs$ which have at least three neighbors in $\doms$ and let $G_1$ be the subgraph of $G'$ formed by deleting all edges between $\subs''$ and $\doms_{< D}$.
    By construction, every vertex in $\subs$ (and hence in $\subs''$) has at most two neighbors in $\doms_{\geq D}$, so $G_1$ has the property that every vertex in $\subs$ has at most two neighbors in $\doms$.
    That is to say, $(G_1,\doms,\subs)$ is a tumor graph.
    By reasoning similar to that behind \cref{eqn:smallbigedges}, we have
    \begin{equation}\label{eqn:smallmededges}
        \abs{\subs''}\leq 2\abs \doms\leq 6n/d\qquad\implies\qquad e\bigl(G'[\subs'',\doms_{< D}]\bigr)\leq 2\abs{\subs''}+2\abs{\doms_{< D}}\leq O(n/d).
    \end{equation}
    \begin{claim}
        With $G_1$ defined as above, $\numb(G',C_{2m+1})\leq\numb(G_1,C_{2m+1})+O(n^{m-1/5})$.
    \end{claim}
    \begin{proof}
        We simply need to bound the number of copies of $C_{2m+1}$ in $G'$ which use some edge within $E\bigl(G[\subs'',\doms_{< D}]\bigr)$.
        In other words, we need to bound the number of copies of $C_{2m+1}$ of the form $(v_1,\dots,v_{2m+1})$ where $v_1\in \doms_{< D}$ and $v_2\in \subs''$.
        In order to bound these, we consider cases according to the nature of $v_4$. \\

        \noindent\textbf{Case: $v_3\in \subs$.}

        Here \cref{eqn:smallmededges} tells us that there are at most $O(n/d)$  choices for the pair $(v_1,v_2)$.
        Then there are at most $\deg v_1<D$ choices for $v_{2m+1}$ and at most $\deg v_2< d$ choices for $v_3$.
        After picking $v_3\in \subs$, we then have at most $\deg v_3< d$ choices for $v_4$.
        Together, this yields a total of at most $O(dDn)$  choices for the tuple $(v_{2m+1},v_1,v_2,v_3,v_4)$.
        Finally, by \Cref{pathasymptotics}, there are at most $2\numb(G',P_{2(m-2)})\leq O(n^{m-2})$  choices for the path $(v_5,\dots,v_{2m})$.
        We conclude that there are at most $O(dDn^{m-1})=O(n^{3/5}\cdot n^{m-1})=O(n^{m-2/5})$  copies of $C_{2m+1}$ of this form.\\

        \noindent\textbf{Case: $v_3\in \doms_{\geq D}$.}

        Here \cref{eqn:smallmededges} tells us that there are at most $O(n/d)$  choices for the pair $(v_1,v_2)$.
        Then there are at most $2$ choices for $v_3$ since $v_2\in \subs$ and $G_1$ is a subgraph of $G'$.
        Therefore, there are at most $O(n/d)$ choices for the triple $(v_1,v_2,v_3)$.
        Finally, by \Cref{pathasymptotics}, there are at most $2\numb(G',P_{2(m-1)})\leq O(n^{m-1})$  choices for the path $(v_4,\dots,v_{2m+1})$.
        We conclude that there are at most $O(n^m/d)=O(n^{m-1/5})$  copies of $C_{2m+1}$ of this form.\\

        \noindent\textbf{Case: $v_3\in \doms_{< D}$.}

        In this situation, \Cref{forks} implies that there are at most $O(\abs{\subs''}+d\abs{\doms_{<D}})\leq O(n)$ many choices for the triple $(v_1,v_2,v_3)$.
        Then there are at most $\deg v_1<D$  choices for $v_{2m+1}$ and at most $\deg v_3< D$  choices for $v_4$.
        This yields a total of at most $O(D^2n)$  choices for the tuple $(v_{2m+1},v_1,v_2,v_3,v_4)$.
        Finally, by \Cref{pathasymptotics}, there are at most $2\numb(G',P_{2(m-2)})\leq O(n^{m-2})$  choices for the path $(v_5,\dots,v_{2m})$.
        We conclude that there are at most $O(D^2n^{m-1})=O(n^{4/5}\cdot n^{m-1})=O(n^{m-1/5})$  copies of $C_{2m+1}$ of this form.
    \end{proof}

    This completes the proof of \Cref{tumorsamecycles}. \hfill\qedsymbol

\subsection{Proof of \Cref{fewbadcycles}: Most cycles in a tumor graph are good}
\label{sec:fewbadcycles}

    Fix a graph $G$ and subsets $V_1,\dots,V_\ell\subseteq V$.
    For $k\geq\ell$, we say that a $k$-cycle $C$ in $G$ contains the pattern $V_1V_2\cdots V_\ell$ if we can cyclically label the vertices of $C$ as $(v_1,\dots,v_k)$ so that $v_i\in V_i$ for all $i\in[\ell]$.
    \medskip

    In order to prove that most copies of $C_{2m+1}$ in a planar tumor graph $(G;\doms,\subs)$ are good, we identify the patterns that a bad copy of $C_{2m+1}$ must contain.
    To do so, we define the following sets:
    \begin{itemize}
        \item $\tp{SSS}$ is the set of all copies of $C_{2m+1}$ that contain the pattern $\subs\subs\subs$.
        \item $\tp{2SS}$ is the set of all copies of $C_{2m+1}$ that contain the pattern $\subs\subs$ at least twice. Note that $\tp{2SS}\supseteq\tp{SSS}$.
        \item $\tp{BBB}$ is the set of all copies of $C_{2m+1}$ containing the pattern $\doms\doms\doms$.
        \item $\tp{BBSS}$ is the set of all copies of $C_{2m+1}$ containing the pattern $\doms\doms\subs\subs$.
        \item $\tp{1BB1SS}$ is the set of all copies of $C_{2m+1}$ containing both the pattern $\doms\doms$ and the pattern $\subs\subs$.
        \item $\tp{2BB}$ is the set of all copies of $C_{2m+1}$ containing the pattern $\doms\doms$ at least twice. Note that $\tp{2BB}\supseteq\tp{BBB}$.
    \end{itemize}

    We note that all bad cycles must fall into at least one of these categories because good cycles contain exactly one instance of either $\doms\doms$ or $\subs\subs$, and otherwise alternate between $\doms$ and $\subs$.
    We now show that each of these six sets is small.

    \begin{lemma}
        The following bounds hold:
        \begin{itemize}
            \item $\abs{\tp{SSS}} \leq O(d^3\abs\subs n^{m-2})$
            \item $\abs{\tp{2SS}\setminus\tp{SSS}} \leq O(d^3\abs\subs n^{m-2})$
            \item $\abs{\tp{BBB}} \leq O(\abs\doms n^{m-1})$
            \item $\abs{\tp{BBSS}\setminus\tp{SSS}} \leq O(\abs\doms n^{m-1})$
            \item $\abs{\tp{1BB1SS}\setminus(\tp{2SS}\cup\tp{BBB}\cup\tp{BBSS})} \leq O(\abs\subs^3 n^{m-4})$
            \item $\abs{\tp{2BB}\setminus(\tp{2SS}\cup\tp{BBB}\cup\tp{BBSS})} \leq O(\abs\subs^3 n^{m-4})$
        \end{itemize}
    \end{lemma}
    \begin{proof}~\\
        \noindent\textbf{Case: \tp{SSS}.} \\
        We need to count cycles of the form $(v_1,\dots,v_{2m+1})$ with $v_2,v_3,v_4\in \subs$.
        To begin, there are at most $2e(G[\subs]) \leq O(\abs\subs)$  choices for the pair $(v_2,v_3)$ after which there are at most $\deg v_2\leq d$  choices for $v_1$ and at most $\deg v_3\leq d$  choices for $v_4$.
        Then, after selecting $v_4$, there are at most $\deg v_4\leq d$  choices for $v_5$.
        Therefore, there are at most $O(d^3\abs\subs)$  choices for the tuple $(v_1,\dots,v_5)$.

        Finally, by \Cref{pathasymptotics}, there are at most $2\numb(G,P_{2(m-2)}) \leq O(n^{m-2})$  choices for the path $(v_6,\dots,v_{2m+1})$.
        In conclusion, there are at most $O(d^3\abs\subs n^{m-2})$  copies of $C_{2m+1}$ of this form. \\

        \noindent\textbf{Case: $\tp{2SS}\setminus\tp{SSS}$.} \\
        We look for a pair of $\subs\subs$ and there will be two paths between them on the cycle, one of odd length and one with even length. Thus, we may assume that the $\subs\subs$ pairs occur at $(v_2,v_3)$ and at $(v_{2i},v_{2i+1})$ for some $i\in\{3,\ldots,m\}$.
        So, we need to count cycles of the form $(v_1,\dots,v_{2m+1})$ with $v_2,v_3,v_{2i},v_{2i+1}\in\subs$.
        There are at most $2e(G[\subs]) \leq O(\abs\subs)$  choices for $(v_2,v_3)$ and for $(v_{2i},v_{2i+1})$.
        Once these are selected, there are at most $d$  choices for each of $v_1,v_4,v_{2i-1}$, which are distinct.
        Therefore, there are at most $O(d^3\abs\subs^2)$  choices for $(v_1,v_2,v_3,v_4,v_{2i-1},v_{2i},v_{2i+1})$.

        Finally, by \Cref{pathasymptotics}, there are at most $2\numb(G,P_{2(i-3)}) \leq O(n^{i-3})$  choices for the path $(v_5,\dots,v_{2i-2})$ and at most $2\numb(G,P_{2(m-i)}) \leq O(n^{m-i})$  choices for the path $(v_{2i+2},\dots,2m+1)$.
        In conclusion, there are at most
        \[
            O\bigl(d^3\abs\subs^2\cdot n^{i-3}\cdot n^{m-i}\bigr) = O\bigl(d^3\abs\subs^2 n^{m-3}\bigr) \leq O\bigl(d^3\abs\subs n^{m-2}\bigr)
        \]
         cycles of this form. \\

        \noindent\textbf{Case: $\tp{BBB}$.} \\
        By \Cref{pathasymptotics}, the number of cycles that contains no $\subs$ vertices at all is at most
        \[
            \numb\bigl(G[\doms],C_{2m+1}\bigr) \leq \numb_\plan(\abs\doms,C_{2m+1}) \leq O(\abs\doms^m) \leq O(\abs \doms n^{m-1})  .
        \]
        Otherwise, we need to count cycles of the form $(v_1,\dots,v_{2m+1})$ with $v_1,v_2,v_3\in\doms$ and $v_4\in\subs$.
        There are at most $2e(G[\doms]) \leq O(\abs\doms)$ choices for $(v_1,v_2)$ and at most $2e(G)\leq O(n)$ choices for $(v_4,v_5)$. Once these are selected, there are at most two choices for $v_3$ because by the definition of a tumor graph, $v_4$ can be adjacent to at most $2$ members of $\doms$.
        Therefore, there are at most $O(\abs\doms n)$ many choices for the tuple $(v_1,\dots,v_5)$.

        Finally, by \Cref{pathasymptotics}, there are at most $2\numb(G,P_{2(m-2)}) \leq O(n^{m-2})$  choices for the path $(v_6,\dots,v_{2m+1})$. In conclusion, there are at most $O(\abs \doms n^{m-1})$ cycles of this form. \\

        \noindent\textbf{Case: $\tp{BBSS}\setminus\tp{SSS}$.} \\
        We need to count cycles of the form $(v_1,\dots,v_{2m+1})$ with $v_1,v_2,v_5\in\doms$ and $v_3,v_4\in\subs$. There are at most $2e(G[\doms]) \leq O(\abs\doms)$ choices for $(v_1,v_2)$ and at most $2e(G[\subs])\leq O(n)$ choices for $(v_3,v_4)$. Once these are selected, there are at most two choices for $v_5$ because by the definition of a tumor graph, $v_4$ can be adjacent to at most $2$ members of $\doms$.

        Finally, by \Cref{pathasymptotics}, there are at most $2\numb(G,P_{2(m-2)}) \leq O(n^{m-2})$  choices for the path $(v_6,\dots,v_{2m+1})$. In conclusion, there are at most $O(\abs \doms n^{m-1})$ cycles of this form. \\

        \noindent\textbf{Case: $\tp{1BB1SS}\setminus(\tp{2SS}\cup\tp{BBB}\cup\tp{BBSS})$.}

        In this case, $m\geq 3$ because any five-cycle containing both the patterns $\doms\doms$ and $\subs\subs$ vertices must be of the form $\doms\doms\subs\subs V$, i.e.\ it belongs to $\tp{BBSS}$.

        Now, any cycle within $\tp{1BB1SS}\setminus(\tp{2SS}\cup\tp{BBB}\cup\tp{BBSS})$ must have the form
        \[
            \doms\subs\subs\doms\subs\ \doms^{k_1}\subs\doms^{k_2}\subs\cdots\doms^{k_\ell}\subs,
        \]
        where $k_1,\dots,k_\ell\in\{1,2\}$ and $\sum_{i=1}^\ell(k_i+1)=2m-4$ and there is at least one $i$ for which $k_i=2$. In fact, for parity reasons, there are at least two $i$'s for which $k_i=2$.

        Therefore, let $t\in[\ell-1]$ be the smallest index for which $k_t=2$.

        If $t=1$, then our cycle has the form
        \[
            \underline{\doms\subs\subs}\ \underline{\doms\subs\doms}\ \underline{\doms\subs\doms}\ \underbrace{V\cdots V}_{2(m-4)}.
        \]
        That is, the initial segment of the cycle $(v_1,v_2,\ldots,v_{2m+1})$ has $v_1,v_4,v_6,v_7,v_9\in\doms$ and has $v_2,v_3,v_5,v_8\in\subs$, and what follows is a path of length $2(m-4)$ with vertices from anywhere in the graph.

        There are at most $4e(G[\subs])$ choices for the initial $\doms\subs\subs$ piece because there are $2e(G[\subs])$ choices for the pair $(v_2,v_3)$ and each vertex in $S$ is adjacent to at most $2$ vertices in $\doms$.
        Similarly, there are at most $2\abs \subs$ choices for each of the $\doms\subs\doms$ pieces. Thus, there are $O(\abs \subs^3)$ many choices for the first 9 vertices.

        Finally, by \Cref{pathasymptotics}, there are at most $2\numb(G,P_{2(m-4)}) \leq O(n^{m-4})$ choices for the remaining vertices.
        Thus there are at most $O(\abs\subs^3n^{m-4})$
        cycles of this form.

        If $t\geq 2$, then our cycle has the form
        \[
            \underline{\doms\subs\subs\doms}\ \underbrace{V\cdots V}_{2(t-2)+1}\ \underline{\doms\subs\doms}\ \underline{\doms\subs\doms}\ \underbrace{V\cdots V\subs}_{2(m-t-3)}.
        \]
        As above, there are at most $8e(G[\subs])\leq O(\abs \subs)$ choices for the $\doms\subs\subs\doms$ piece and at most $2\abs\subs$ choices for each of the $\doms\subs\doms$ pieces.

        Finally, by \Cref{pathasymptotics}, there are at most $2\numb(G,P_{2(t-2)+1)}) \leq O(n^{t-1})$ choices for the path between the $\doms\subs\subs\doms$ piece and the first $\doms\subs\doms$ piece and there are at most $2\numb(G,P_{2(m-t-3)}) \leq O(n^{m-t-3})$ choices for the path following the second $\doms\subs\doms$ piece.
        Thus there are at most $O(\abs\subs^3 n^{m-4})$
        cycles of this form for any integer $t\in\{1,\ldots,\ell-1\}$. \\

        \noindent\textbf{Case: $\tp{2BB}\setminus(\tp{2SS}\cup\tp{BBB}\cup\tp{BBSS})$.}

        In this case, $m\geq 3$ because any five-cycle with two instances of $\doms\doms$ must be of the form $\doms\doms\doms$.
        Any cycle within $\tp{2BB}\setminus(\tp{2SS}\cup\tp{BBB}\cup\tp{BBSS})$ must have the form
        \[
            \doms\subs\doms^{k_1}\subs\doms^{k_2}\subs\cdots \doms^{k_\ell}\subs\doms,
        \]
        where $k_1,\dots,k_\ell\in\{1,2\}$ and $\sum_{i=1}^\ell(k_i+1)=2m-2\geq 4$ and there is an $i$ for which $k_i=2$. In fact, for parity reasons, there are at least two $i$'s for which $k_i=2$.

        If $k_i=2$ for all $i\in\{1,\ldots,\ell\}$, then $3\ell=2m+1$ (hence $m\geq 4$) and our cycle has the form
        \[
            \underbrace{\underline{\doms\subs\doms}\ \underline{\doms\subs\doms}\ \cdots\ \underline{\doms\subs\doms}}_\ell.
        \]
        There are at most $2\abs \subs$ choices for each $\doms\subs\doms$ piece.
        Hence there are at most
        \[
            O\bigl(\abs \subs^{(2m+1)/3}\bigr) \leq O\bigl(\abs \subs^3 n^{(2m-8)/3}\bigr) \leq O(\abs \subs^3 n^{m-4})
        \]
        cycles of this form for all $m\geq 4$.

        If not all $k_i$'s are equal to $2$, we may assume, without loss of generality, that $k_{\ell}=1$ and there exists a $t\in\{1,\ldots,\ell-2\}$ which is the smallest index for which $k_t=2$.

        If $t=1$, then our cycle has the form
        \[
            \underline{\doms\subs\doms}\ \underline{\doms\subs\doms}\ \underbrace{V\cdots V\subs}_{2(m-4)}\ \underline{\doms\subs\doms}.
        \]
        There are at most $2\abs \subs$ choices for each of the $\doms\subs\doms$ pieces.

        Finally, by \Cref{pathasymptotics}, there are at most $2\numb(G,P_{2(m-4)}) \leq O(n^{m-4})$ choices for the remaining vertices.
        Thus there are at most $O(\abs\subs^3 n^{m-4})$
        cycles of this form.

        If $t\geq 2$, then our cycle has the form
        \[
            \underline{\doms\subs\doms}\ \underbrace{\subs\doms\cdots \subs\doms}_{2(t-1)}\ \underline{\doms\subs\doms}\ \underbrace{V\cdots V\subs}_{2(m-t-3)}\ \underline{\doms\subs\doms}.
        \]
        As above, there are at most $2\abs \subs$ choices for each of the $\doms\subs\doms$ pieces.

        Finally, by \Cref{pathasymptotics}, there are at most $2\numb(G,P_{2(t-1)}) \leq O(n^{t-1})$ choices for the path between the first $\doms\subs\doms$ piece and the second $\doms\subs\doms$ piece and there are at most $2\numb(G,P_{2(m-t-3)}) \leq O(n^{m-t-3})$ choices for the path between the second $\doms\subs\doms$ piece and the last $\doms\subs\doms$ piece.
        Thus there are at most $O(\abs\subs^3 n^{m-4})$ cycles of this form for any integer $t\in\{1,\ldots,\ell-2\}$. \\
    \end{proof}

    This completes the proof of \Cref{fewbadcycles}. \hfill\qedsymbol

\subsection{Proof of \Cref{cleaning}: Cleaning a tumor graph}
\label{sec:cleaning}

The process of cleaning $G_1$ to arrive at $G_2$ will go through several stages and requires a number of facts.

In \Cref{tumor:fewcrossedges}, we establish that there are few edges between distinct tumors. In \Cref{tumor:conuncon}, we establish that deleting edges between vertices in $\subs$ does not decrease the number of good cycles by much and establish that an operation we call ``contraction--uncontraction'' maintains planarity and the number of vertices but does not decrease the number of good cycles at all.

In \Cref{tumor:onemor} (Cleaning Stage I), we delete certain edges and vertices that cannot participate in good cycles and perform contraction--uncontraction on some edges between vertices in $\bigsqcup_{x\in\doms}\onemor x$. In \Cref{tumor:onemortumor} (Cleaning Stage II), we delete and contract--uncontract some edges such that every vertex in $\bigsqcup_{x\in\doms}\onemor x$ has at most one neighbor in a tumor and there are no edges between distinct tumors. In \Cref{tumor:tumors} (Cleaning Stage III), we perform deletion and contraction--uncontraction and modify $\subs$ and $\doms$ into $\subs'$ and $\doms'$ so that any edge induced by $\subs'$ is within some tumor. In each stage, we ensure that the total number of good cycles does not change by too much.

\subsubsection{Few edges between distinct tumors}
\label{tumor:fewcrossedges}

We begin with two simple observations about planar tumor graphs.

\begin{prop}\label{few_in_tumor}
    Let $(G;\doms,\subs)$ be a planar tumor graph.
    \begin{itemize}
        \item For any $xy\in\binom{\doms}{2}$ and any $z\in V(G)\setminus\{x,y\}$, we have $\abs{N(z)\cap\tumor xy}\leq 2$.
            In particular, $e\bigl(G[\subs_{xy}])\leq \abs{\subs_{xy}}$.
        \item $G$ has at most four edges between any fixed pair of distinct tumors.
    \end{itemize}
\end{prop}
\begin{proof}
    If $\abs{N(z)\cap\tumor xy}\geq 3$, then $x,y,z$ would be three distinct vertices with three common neighbors,  so $G$ would contain a copy of $K_{3,3}$, contradicting the fact that $G$ is planar.
    This establishes the first item.
    \medskip

    To prove the second item, fix two distinct tumors $\tumor xy,\tumor zw$ of $G$ and set $T=G[\tumor xy,\tumor zw]$.
    Since $xy\neq zw$, we may relabel these vertices so that $x,y,z$ are distinct.
    Now, suppose that $T$ contained at least five edges; by the first item, the maximum degree of $T$ is at most two and so $T$ must contain a matching on three edges.
    Let $a_1b_1,a_2b_2,a_3b_3\in E(T)$ be such a matching where $a_i\in\tumor xy$ and $b_i\in\tumor zw$ for each $i\in[3]$.
    But then these six vertices along with $x,y,z$ contain a subdivided copy of $K_{3,3}$ with parts $\{x,y,z\}$ and $\{a_1,a_2,a_3\}$; a contradiction.
\end{proof}

The goal of the remainder of this section is to prove that planar tumor graphs have few edges between distinct tumors (\Cref{tumor-tumor}) and that most other vertices in $\subs$ interact sparsely with the set of tumors (\Cref{onemor-tumor})

\begin{lemma}\label{tumor-tumor}
    If $(G;\doms,\subs)$ is a planar tumor graph, then $G$ has at most $O(\abs\doms)$ edges of the form $uv$, with $u$ and $v$ in distinct tumors.
\end{lemma}

\begin{lemma}\label{onemor-tumor}
    Let $(G;\doms,\subs)$ be a planar tumor graph.
    Define $X$ to be the set of all vertices $v$ with any of the following properties:
    \begin{enumerate}
        \item $v\in\nomor\sqcup\bigsqcup_{x\in\doms}\onemor x$ and $v$ has at least three neighbors within $\bigsqcup_{xy\in{\doms\choose 2}}\tumor xy$.
        \item $v\in\bigsqcup_{x\in\doms}\onemor x$ and $v$ has neighbors in distinct tumors.
        \item $v\in\onemor x$ for some $x\in\doms$ and $v$ has two neighbors within $\tumor yz$ for some $yz\not\ni x$.
    \end{enumerate}
    Then
    \begin{itemize}
        \item $\abs X \leq O(\abs \doms)$, and
        \item The number of edges between $X$ and $\bigsqcup_{xy\in{\doms\choose 2}}\tumor xy$ is at most $O(\abs\doms)$.
    \end{itemize}
\end{lemma}

In order to prove \Cref{tumor-tumor,onemor-tumor}, we need the notion of a separation of a tumor graph.

To a tumor graph $(G;\doms,\subs)$, we associate an auxiliary graph $\tum(G;\doms,\subs)$ which has vertex set $\doms$ and $xy$ is an edge whenever $\subs_{xy}\neq\varnothing$.

\begin{defn}
    A tumor graph $(G;\doms,\subs)$ is said to be \emph{separated} if $\tum(G;\doms,\subs)$ is a matching, possibly with isolated vertices.
\end{defn}

Unsurprisingly, separated tumor graphs are much easier to handle and so we will need to ``separate'' a tumor graph in order to obtain the bounds in the preceding lemmas.

\begin{defn}\label{defseparation}
    Let $(G;\doms,\subs)$ be a tumor graph.
    A tumor graph $(G';\doms',\subs)$ is said to be a \emph{separation} of $(G;\doms,\subs)$ if
    \begin{itemize}
        \item $(G';\doms',\subs)$ is separated, and
        \item $G[\subs]=G'[\subs]$, and
        \item There is a partition $\doms'=\bigsqcup_{x\in\doms}\doms'_x$ with the following properties for any $s\in\subs$ and any $xy\in{\doms\choose 2}$:
            \begin{itemize}
                \item If $s\in\onemor x$, then $s\in\onemor{x'}$ for some $x'\in\doms'_x$, and
                \item If $s\in\tumor xy$, then $s\in\tumor{x'}{y'}$ for some $x'\in\doms'_x$ and $y'\in\doms'_y$.
            \end{itemize}
    \end{itemize}
\end{defn}

That is, each $x\in\doms$ corresponds to a set of vertices $\doms'_x\subset\doms'$, and the vertices of $\subs$ adjacent to $x$ in $G$ must be adjacent to a unique member of $\doms'_x$ in $G'$.

A separation of a tumor graph preserves the following crucial property:
\begin{obs}\label{obssep}
    Let $(G;\doms,\subs)$ be a tumor graph and let $(G';\doms',\subs)$ be a separation.
    If $u,v\in\subs$ are in distinct clusters\footnote{A ``cluster'' here refers to any subset of $\subs$ of the form $\tumor xy$, $\onemor x$ or $\nomor$.} of $G$, then $u,v$ are in distinct clusters of $G'$.
\end{obs}
Seeing how both \Cref{tumor-tumor,onemor-tumor} bound interactions between distinct clusters, the first step in their proofs will be to find a separation of the original tumor graph that is not too large.
This is the content of the following proposition.

\begin{prop}[Separation Proposition]\label{separation}
    Any planar tumor graph $(G;\doms,\subs)$ has a planar separation $(G';\doms',\subs)$ where
    \[
        \abs{\doms'}\leq 2e\bigl(\tum(G;\doms,\subs)\bigr)+\abs{\{x\in\doms:\deg_{\tum(G;\doms,\subs)} x = 0\}}\leq 6\abs\doms.
    \]
\end{prop}
\begin{proof}
    Set $\tum=\tum(G;\doms,\subs)$
    Note that $G$ contains a subdivision of $\tum$ and so $\tum$ is additionally planar.
    Due to this, the bound of $2e(\tum)+\abs{\{x\in\doms:\deg_{\tum} x = 0\}}\leq 6\abs\doms$ is immediate, so we focus only on the first inequality

    We prove the claim by double induction on the pair $(\Delta,\eta)$ (induction is done on $\Delta$ first and then $\eta$) where $\Delta$ is the maximum degree of $\tum$ and $\eta$ is the number of vertices in $\tum$ of degree $\Delta$.

    If $\Delta\leq 1$, then $(G;\doms,\subs)$ is separated and so there is nothing to prove.
    Thus, suppose that $\Delta\geq 2$.

    To begin, we may suppose that $G$ has no edges between vertices of $\doms$ since we may remove any such edges without affecting the conclusion of the lemma.
    Now, fix a straight-edge planar embedding of $G$ and let $x\in\doms$ be any vertex with $\deg_\tum x=\Delta$.
    We may label the neighbors of $x$ as $s_0,\dots,s_{k-1}$ in counter-clockwise order around $x$.
    Since $G$ has no edges between vertices of $\doms$, each $s_i$ resides within $\subs$.
    Let $\{t_0,\dots,t_{\ell-1}\}=\{s_0,\dots,s_{k-1}\}\setminus\onemor x$ where the $t_i$'s remain in counter-clockwise order about $x$; note that each $t_i$ resides within $\tumor xy$ for some $y\in\doms$.
    Note that $\ell\geq 2$ since $x$ has degree $\Delta\geq 2$ in $\tum$.

    We claim that there is some $y\in\doms$ for which $\{t_0,\dots,t_{\ell-1}\}\cap\tumor xy$ is a non-trivial cyclic interval; that is $\{t_0,\dots,t_{\ell-1}\}\cap\tumor xy=\{t_i,t_{i+1},\dots,t_{i+r}\}$ for some $i\in\{0,\dots,\ell-1\}$ and $r\in[\ell-1]$ where the indices are computed modulo $\ell$.
    If no such $y$ were to exist, then we could locate indices $a<b<c<d\in\{0,\dots,\ell-1\}$ for which $t_a,t_c\in\tumor xy$ and $t_b,t_d\in\tumor xz$ for some distinct $y,z\in\doms\setminus\{x\}$.
    Now, consider the subgraph of $G$ induced by $x,y,z,t_a,t_b,t_c,t_d$, inheriting its planar embedding from $G$; call this plane graph $H$.
    In $H$, the neighbors of $x$ are $t_a,t_b,t_c,t_d$, which appear in counter-clockwise order around $x$. By a standard argument in planar graph theory, we may suppose that $t_at_bt_ct_d$ forms a cycle in $H$ since we can add these edges without violating the planarity of $H$.
    However, as in \Cref{fig:counterclockwise}, $t_ay,yt_c$ and $t_bz,zt_d$ are edges of $H$, and so $H$ is a subdivision of $K_5$; a contradiction.

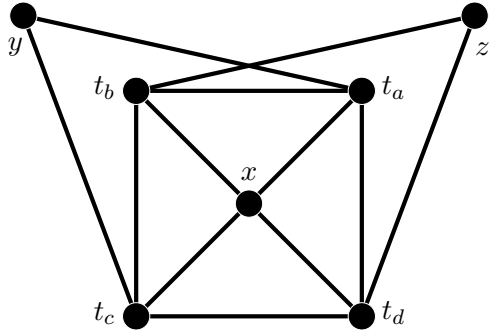
\begin{figure}[ht]
    \centering
        \begin{tikzpicture}
            \useasboundingbox (-3.3,-1.7) rectangle (3.3,2.7);
            \coordinate (x) at (0,0);
            \coordinate (y) at (-3,2.5);
            \coordinate (z) at (3,2.5);
            \coordinate (ta) at (1.5,1.5);
            \coordinate (tb) at (-1.5,1.5);
            \coordinate (tc) at (-1.5,-1.5);
            \coordinate (td) at (1.5,-1.5);
            \begin{pgfonlayer}{fore}
                \node[vertex,label={[yshift=0pt]$x$}] at (x) {};
                \node[vertex,label={[xshift=-3pt, yshift=-24pt]$y$}] at (y) {};
                \node[vertex,label={[xshift=3pt, yshift=-24pt]$z$}] at (z) {};
                \node[vertex,label={[xshift=12pt, yshift=-12pt]$t_a$}] at (ta) {};
                \node[vertex,label={[xshift=-12pt, yshift=-12pt]$t_b$}] at (tb) {};
                \node[vertex,label={[xshift=-12pt, yshift=-12pt]$t_c$}] at (tc) {};
                \node[vertex,label={[xshift=12pt, yshift=-12pt]$t_d$}] at (td) {};
            \end{pgfonlayer}
            \begin{pgfonlayer}{back}
                \foreach\i in {a,b,c,d}
                    \draw[edge] (x) -- (t\i);
                \draw[edge] (ta) -- (tb) -- (tc) -- (td) -- (ta);
                \draw[edge] (ta) -- (y) -- (tc);
                \draw[edge] (tb) -- (z) -- (td);
            \end{pgfonlayer}
        \end{tikzpicture}
        \caption{\label{fig:counterclockwise}The vertices $t_a,t_c\in\subs_{xy}$ and $t_b,t_c\in\subs_{xz}$ in a counterclockwise order. The vertices $\{x,t_a,t_b,t_c,t_d\}$ form the vertices of a subdivision of $K_5$.}
\end{figure}

    Thus, without loss of generality, let $y\in\doms$ be such that $\{t_0,\dots,t_{\ell-1}\}\cap\tumor xy=\{t_0,\dots,t_r\}$ for some $r\in[\ell-1]$.
    We may additionally suppose that $t_0=s_0$ and that $t_r=s_{r'}$ for some $r'\in[k-1]$.
    We form the new tumor graph $(G';\doms',\subs)$ by introducing a new vertex $x'$ to have $\doms'=\doms\sqcup\{x'\}$ and replacing all edges of the form $s_ix$ by $s_ix'$ for each $i\in\{0,\dots,r'\}$ (see \Cref{fig:separation}).
    Observe that $G'$ is still planar since $s_0,\dots,s_{r'}$ is a cyclic interval of neighbors of $x$.

\newcommand{\separateplot}[1]{%
        \begin{tikzpicture}
            \def\rad{2}
            \pgfmathsetmacro{\ycor}{#1*\rad};
            \coordinate (xp) at (0,\ycor);
            \pgfmathsetmacro{\test}{int(1000*\ycor}
            \begin{pgfonlayer}{fore}
                \ifnum\test=0
                \else
                    \node[vertex,label={[xshift=12pt,yshift=-12pt]$x'$}] at (xp) {};
                \fi
            \end{pgfonlayer}
            \useasboundingbox (-1.7*\rad,-1.5*\rad) rectangle (1.7*\rad,3.1*\rad);
            \coordinate (x) at (0,0);
            \coordinate (y) at (0,\ycor+2.5*\rad);
            \foreach \i/\th in {0/40, 1/85, 2/115, 3/130, 4/140}
                \coordinate (s\i) at ({\rad*tan(90-\th)},{\ycor+1.25*\rad});
            \foreach \i/\th in {1/50, 2/68, 4/100}
                \coordinate (t\i) at ({0.65*\rad*tan(90-\th)},{\ycor+0.65*1.25*\rad});
            \foreach \i/\th in {0/5, 1/25}
                \coordinate (e\i) at ({0.9*\rad*cos(\th)/cos(15-\th)},{0.9*\rad*sin(\th)/cos(15-\th)});
            \foreach\i/\th in {0/150, 1/165, 2/180}
                \coordinate (f\i) at ({0.9*\rad*cos(\th)/cos(165-\th)},{0.9*\rad*sin(\th)/cos(165-\th)});
            \coordinate (za) at (320:2*\rad);
            \foreach \i/\th in {0/290, 1/330, 2/350}
                \coordinate (a\i) at ({\rad*cos(\th)/cos(320-\th)},{\rad*sin(\th)/cos(320-\th)});
            \foreach \i/\th in {1/310}
                \coordinate (c1) at ({0.75*\rad*cos(\th)/cos(320-\th)},{0.75*\rad*sin(\th)/cos(320-\th)});
            \coordinate (zb) at (220:2*\rad);
            \foreach \i/\th in {0/190, 1/210, 2/230, 3/250}
                \coordinate (b\i) at ({\rad*cos(\th)/cos(220-\th)},{\rad*sin(\th)/cos(220-\th)});
            \foreach \i/\th in {0/200, 1/240}
                \coordinate (d\i) at ({0.75*\rad*cos(\th)/cos(220-\th)},{0.75*\rad*sin(\th)/cos(220-\th)});
            \begin{pgfonlayer}{fore}
                \node[vertex,label={[yshift=-27pt]$x$}] at (x) {};
                \node[vertex,label={[xshift=12pt,yshift=-12pt]$y$}] at (y) {};
                \foreach\i in {0,1,2,3,4}
                    \node[vertex] at (s\i) {};
                \node[label={[xshift=12pt,yshift=0pt]$s_0\!=\!t_0$}] at (s0) {};
                \node[label={[xshift=17pt,yshift=0pt]$s_1\!=\!t_3$}] at (s1) {};
                \node[label={[xshift=-10pt,yshift=0pt]$s_4\!=\!t_7$}] at (s4) {};
                \foreach\i in {1,2,4}
                    \node[label={[xshift=0pt,yshift=0pt]$t_\i$}] at (t\i) {};
                \foreach\i in {1,2,4}
                    \node[vertex] at (t\i) {};
                \node[vertex] at (za) {};
                \foreach\i in {0,1,2}
                    \node[vertex] at (a\i) {};
                \node[vertex] at (c1) {};
                 \node[vertex] at (zb) {};
               \foreach\i in {0,1,2,3}
                    \node[vertex] at (b\i) {};
                \foreach\i in {0,1}
                    \node[vertex] at (d\i) {};
                \foreach\i in {0,1}
                    \node[vertex] at (e\i) {};
                \foreach\i in {0,1,2}
                    \node[vertex] at (f\i) {};
            \end{pgfonlayer}
            \begin{pgfonlayer}{back}
                \foreach\i in {0,1,2,3,4}
                    \draw[edge] (xp) -- (s\i) -- (y);
                \foreach\i in {1,2,4}
                    \draw[edge] (xp) -- (t\i);
                \foreach\i in {0,1,2}
                    \draw[edge] (x) -- (a\i) -- (za);
                \draw[edge] (x) -- (c1);
                \foreach\i in {0,1,2,3}
                    \draw[edge] (x) -- (b\i) -- (zb);
                \foreach\i in {0,1}
                    \draw[edge] (x) -- (d\i);
               \foreach\i in {0,1}
                    \draw[edge] (x)--(e\i);
                \foreach\i in {0,1,2}
                    \draw[edge] (x)--(f\i);
                \draw[edge] (s0) -- (s1) -- (s2) -- (s3) -- (s4);
                \draw[edge] (a0) -- (a1) -- (a2);
                \draw[edge] (b0) -- (b1) -- (b2) -- (b3);
            \end{pgfonlayer}
        \end{tikzpicture}
}

    \begin{figure}[ht]
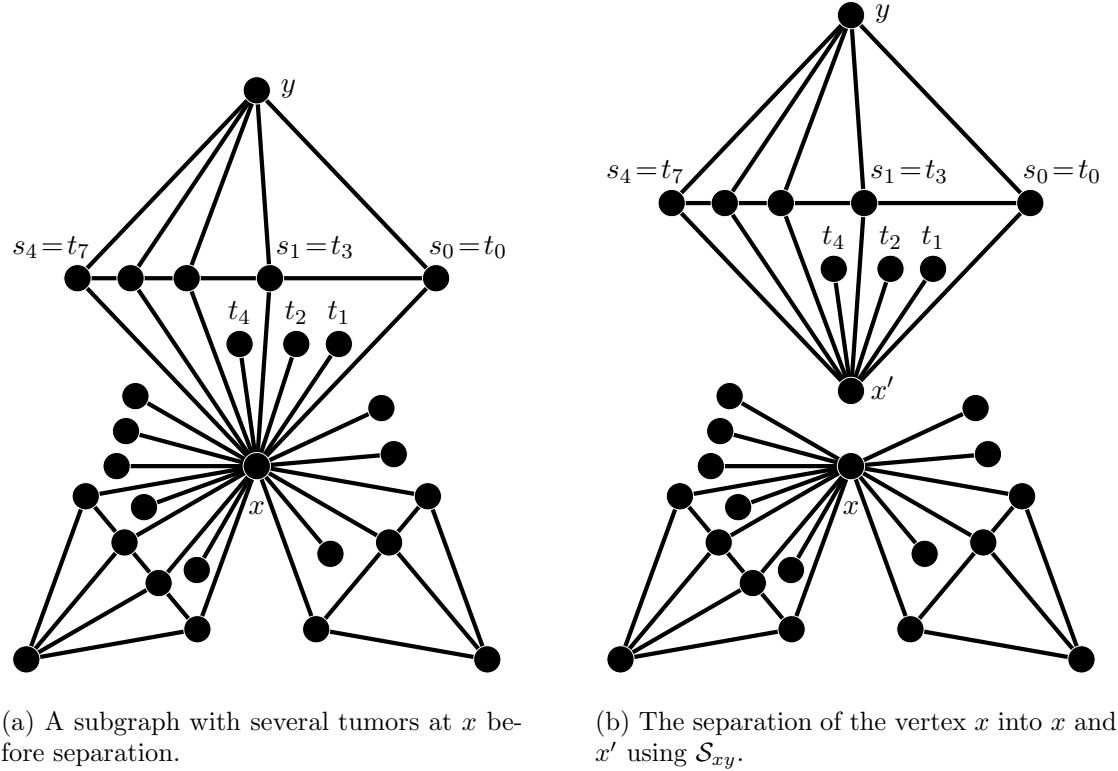

        \centering
        \begin{subfigure}[b]{0.42\textwidth}
            \separateplot{0}
            \caption{\label{fig:separation:before}A subgraph with several tumors at $x$ before separation.}
        \end{subfigure}
        \hfil
        \begin{subfigure}[b]{0.42\textwidth}
            \separateplot{0.5}
            \caption{\label{fig:separation:after}The separation of the vertex $x$ into $x$ and $x'$ using $\subs_{xy}$.}
        \end{subfigure}
        \caption{\label{fig:separation}An example of the splitting of vertices in the proof of \Cref{separation}.}
    \end{figure}

    Now, certainly $(G';\doms',\subs)$ is a tumor graph and $G[\subs]=G'[\subs]$.
    Furthermore, setting $\tum'=\tum(G';\doms',\subs)$, we find that $\deg_{\tum'}x'=1$ and $\deg_{\tum'}x=\deg_{\tum}x-1$.
    In particular,
    \[
        2e(\tum)+\abs{\{x\in\doms:\deg_{\tum} x = 0\}}=2e(\tum')+\abs{\{x\in\doms:\deg_{\tum'} x = 0\}}.
    \]
    Furthermore, if $\Delta'$ denotes the maximum degree of $\tum'$ and $\eta'$ denotes the number of vertices in $\tum'$ of degree $\Delta'$, then we find that $(\Delta',\eta')$ is strictly smaller than $(\Delta,\eta)$ in the lexicographic ordering.
    Thus, the claim follows from the induction hypothesis.
\end{proof}

Now that we understand how to separate a tumor graph, we can prove \Cref{tumor-tumor,onemor-tumor}.
Both proofs follow the same philosophy: separate, contract, bound.

\begin{proof}[Proof of \Cref{tumor-tumor}]
    Let $(G';\doms',\subs)$ be the separation of $(G;\doms,\subs)$ guaranteed by \Cref{separation}.
    This guarantees that edges between distinct tumors of $G$ are between distinct tumors of $G'$ (\Cref{obssep}).
    Moreover, $|\doms'|\leq 6|\doms|$, so it suffices to prove the claim for $(G';\doms',\subs)$.
    In other words, we may suppose that $(G;\doms,\subs)$ is already separated.
    \medskip

    Let $R$ denote the set of edges with end-points in distinct tumors and set $\tum=\tum(G;\doms,\subs)$, which is a matching since $(G;\doms,\subs)$ is separated; in particular, $e(\tum)\leq\abs\doms/2$.

    Now, create a graph $H$ whose vertex set is $E(\tum)$ where $\{xy,zw\}\in E(H)$ if there is an edge between $\tumor xy$ and $\tumor zw$ in $G$.
    Due to \Cref{few_in_tumor}, we know that $\abs R\leq 4e(H)$.
    The key observation is that $H$ is a planar graph.
    Indeed, consider starting with $G$ and contracting all edges of the form $sb\in E(G)$ for $s\in\bigsqcup_{xy\in{\doms\choose 2}}\tumor xy$ and $b\in\doms$.
    Since $\tum$ is a matching, the effect of these contractions is to replace each tumor by a single vertex; in particular, $H$ is isomorphic to a subgraph of this contracted graph.

    Putting these observations together, we finally bound
    \[
        {1\over 4}\abs R\leq e(H)\leq 3v(H)=3e(\tum)\leq {3\over 2}\abs\doms\qquad\implies\qquad \abs R\leq 6\abs\doms.\qedhere
    \]
\end{proof}

The proof of \Cref{onemor-tumor} follows along similar lines, but is more involved since we will need to perform many different sequences of contractions.

\begin{proof}[Proof of \Cref{onemor-tumor}]
    Let $(G';\doms',\subs)$ be the separation of $(G;\doms,\subs)$ guaranteed by \Cref{separation}.
    Since $\abs{\doms'}\leq 6\abs\doms$,  it suffices to prove the claim for $(G';\doms',\subs)$ (see \Cref{obssep}).
    In other words, we may suppose that $(G;\doms,\subs)$ is already separated.

    Begin by fixing any $xy\in{\doms\choose 2}$ and consider $\tumor xy$.
    We build an auxiliary graph $H_{xy}$ whose vertex set is $\tumor xy$ where $ab$ is an edge if there is some $s\in\nomor\sqcup\bigsqcup_{z\in\doms}\onemor z$ for which $s$ is adjacent to both $a$ and $b$ in $G$.
    Observe that any such $s$ is adjacent to at most two vertices in $\tumor xy$ (by \Cref{few_in_tumor}), hence $G$ has a subdivision of $H_{xy}$ and so $H_{xy}$ is a planar graph.
    In particular, the chromatic number of $H_{xy}$ is bounded by some absolute constant $C$.\footnote{The actual value of $C$ is inconsequential to the proof. The four color theorem implies that $C\leq 4$, though the easier bound of $C\leq 6$ would suffice. In fact, one can show that $C\leq 3$ in this special case.}

    We may therefore fix a coloring $\chi\colon\bigsqcup_{xy\in{\doms\choose 2}}\tumor xy\to[C]$ so that $\chi$ is a proper coloring of each $H_{xy}$.
    Next, fix an arbitrary orientation $\overrightarrow\tum$ of the matching $\tum(G;\doms,\subs)$.
    For each $t\in[C]$, we build a graph $G_t$ from $G$ as follows: for each $(x,y)\in\overrightarrow\tum$,
    \begin{itemize}
        \item Contract all edges of the form $xa$ where $a\in\tumor xy$ and $\chi(a)=t$, and
        \item Contract all edges of the form $ya$ where $a\in\tumor xy$ and $\chi(a)\neq t$.
    \end{itemize}

    We will use the graphs $G_t$ to define sets $A_t$ such that $X\subseteq A_1\cup\cdots\cup A_C$ and then show that each $A_t$ has size at most $2\abs \doms$ and that the number of edges between $A_t$ and $\bigsqcup_{xy\in{\doms\choose 2}}\tumor xy$ is at most $6\abs \doms$.

    To that end, for each $t\in[C]$, let $A_t\subseteq \nomor\sqcup\bigsqcup_{x\in\doms}\onemor x$ denote those vertices that have at least three neighbors within $\doms$ in the graph $G_t$. (Note that no vertices of $\nomor\sqcup\bigsqcup_{x\in\doms}\onemor x$ were lost when creating $G_t$.)

    We first consider those $v\in\nomor\sqcup\bigsqcup_{x\in\doms}\onemor x$ which have at least three neighbors within $\bigsqcup_{xy\in{\doms\choose 2}}\tumor xy$.

    Suppose first that $v$ has neighbors within at least three distinct tumors: $\tumor{x_1}{y_1},\tumor{x_2}{y_2},\tumor{x_3}{y_3}$.
    Since $G$ is separated, each of the $x_i$'s and $y_i$'s are distinct.
    As such, in each $G_t$, $v$ is adjacent to either $x_i$ or $y_i$ (or both) for each $i\in[3]$ and so $v\in A_t$ for each $t\in[C]$.

    If this is not the case, then since $v$ has at most two neighbors within any individual $\tumor xy$ (by \Cref{few_in_tumor}), this means that $v$ has neighbors within two distinct tumors, $\tumor{x_1}{y_1},\tumor{x_2}{y_2}$, such that it has two neighbors $a,b\in\tumor{x_1}{y_1}$.
    Again, since $G$ is separated, $x_1,y_1,x_2,y_2$ are distinct.
    Now, since $v$ is a common neighbor of $a$ and $b$, we know that $ab\in E(H_{x_1y_1})$ and so $\chi(a)\neq\chi(b)$.
    In particular, if $t=\chi(a)$, then $vx_1$ and $vy_1$ are both edges of $G_t$.
    Finally, since $v$ has a neighbor in $\tumor{x_2}{y_2}$, either $vx_2$ or $vy_2$ is an edge of $G_t$ and so $v\in A_t$.
    \vskip5pt

    Next, suppose that $v\in\onemor x$ for some $x\in\doms$ and suppose that $v$ has neighbors within distinct tumors $\tumor yz$ and $\tumor wa$.
    Since $G$ is separated, we know that $y,z,w,a$ are distinct; in particular, at most one of these four vertices is equal to $x$.
    Without loss of generality, we may suppose that $x\not\in\{y,z,w\}$.
    Now, since $v$ has some neighbor within $\tumor wa$, there is some value of $t$ for which $vw$ is an edge of $G_t$.
    Within this same $G_t$, either $vy$ or $vz$ is also an edge.
    Finally, $vx$ is additionally an edge of $G_t$ and so $v\in A_t$.
    \vskip5pt

    Finally, suppose that $v\in\onemor x$ for some $x\in\doms$ and that $v$ has two neighbors within $\tumor yz$ for some $\{y,z\}\not\ni x$; call these two neighbors $a,b$.
    As before, we know that $ab\in E(H_{yz})$ and so $\chi(a)\neq \chi(b)$.
    Thus, if $t=\chi(a)$, then both $vy$ and $vz$ are edges of $G_t$.
    Additionally, $vx$ is an edge of $G_t$ and so $v\in A_t$ since $x,y,z$ are distinct. Thus $X\subseteq \bigcup_{t\in [C]} A_t$.
    \vskip5pt

     Since $G_t$ is planar and each vertex in $A_t$ has at least three neighbors in $\doms$ and $A_t$ is disjoint from $\doms$, \Cref{bipartite} implies that $\abs{A_t}\leq 2\abs\doms$.
    Next, since each $v\in A_t$ had at most two neighbors in $G$ within any particular tumor (\Cref{few_in_tumor}), the number of edges between $A_t$ and $\bigsqcup_{xy\in{\doms\choose 2}}\tumor xy$ in $G$ is at most twice as large as the number of edges between $A_t$ and $\doms$ in $G_t$.
    Of course, the number of edges between $A_t$ and $\doms$ in $G_t$ is at most $2(\abs{A_t}+\abs{\doms})\leq 6\abs\doms$.

    Putting these together, we have shown that $|X|\leq C\cdot 2\abs\doms$ and the number of edges between $X$ and $\bigsqcup_{xy\in{\doms\choose 2}}\tumor xy$ is at most $C\cdot 6\abs\doms$.
    Since $C$ is bounded, the claim follows.
\end{proof}

We will not need the notion of separation in the rest of the proof but will instead use \Cref{tumor-tumor,onemor-tumor} to conclude structural facts about the graph $(G;\doms,\subs)$.

\subsubsection{Contraction--uncontraction}
\label{tumor:conuncon}

In this section, we introduce the main operation used to control tumor graphs: contraction--uncontraction.

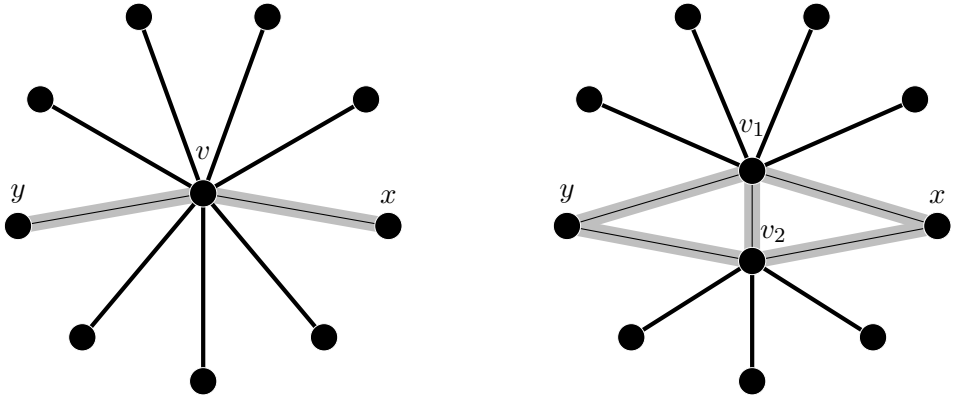
\begin{figure}[ht]
    \begin{center}
        \begin{tikzpicture}
            \def\num{9};
            \def\rad{2.5};
            \def\tilt{10};
            \coordinate (center) at (0,0);
            \foreach\i in {0,...,{\num-1}}
                \coordinate (\i) at ({\i*360/\num-\tilt}:\rad);
            \coordinate (x) at (0) {};
            \coordinate (y) at (5) {};
            \begin{pgfonlayer}{fore}
                \foreach\i in {0,...,{\num-1}}
                    \node[vertex] at (\i) {};
                \node[vertex,label={[yshift=4pt]$v$}] at (center) {};
                \node[label={[yshift=1pt]$x$}] at (x) {};
                \node[label={[yshift=1pt]$y$}] at (y) {};
            \end{pgfonlayer}
            \begin{pgfonlayer}{main}
                \foreach\x in {1,2,3,4,6,7,8}
                    \draw[edge] (center)--(\x);
                \draw [thedge] (x) -- (center) -- (y);
            \end{pgfonlayer}
            \begin{pgfonlayer}{back}
                \draw [spedge] (x) -- (center) -- (y);
            \end{pgfonlayer}
        \end{tikzpicture}\hfil
        \begin{tikzpicture}
            \def\num{9}
            \def\split{5}
            \def\rad{2.5};
            \def\tilt{10};
            \coordinate (center1) at (0,0.3);
            \coordinate (center2) at (0,-0.9);
            \foreach\i in {0,...,\num}
                \coordinate (\i) at ({\i*360/\num-\tilt}:\rad);
            \coordinate (x) at (0) {};
            \coordinate (y) at (5) {};
            \begin{pgfonlayer}{fore}
                \foreach\i in {0,...,\num}
                    \node[vertex] at (\i) {};
                \node[vertex,label={[yshift=4pt]$v_1$}] at (center1) {};
                \node[vertex,label={[xshift=8pt,yshift=-2pt]$v_2$}] at (center2) {};
                \node[label={[yshift=1pt]$x$}] at (x) {};
                \node[label={[yshift=1pt]$y$}] at (y) {};
            \end{pgfonlayer}
            \begin{pgfonlayer}{main}
                \draw[thedge] (x)--(center1)--(y);
                \draw[thedge] (x)--(center2)--(y);
                \draw[thedge] (center1)--(center2);
                \foreach\i in {1,2,3,4,6,7,8}
                    \ifthenelse {\i<\split}
                                {\draw[edge] (center1)--(\i)}
                                {\draw[edge] (center2)--(\i)};
            \end{pgfonlayer}
            \begin{pgfonlayer}{back}
                \draw[spedge] (x)--(center1)--(y);
                \draw[spedge] (x)--(center2)--(y);
                \draw[spedge] (center1)--(center2);
            \end{pgfonlayer}
        \end{tikzpicture}
    \end{center}
\caption{Uncontracting along the path $xvy$.\label{uncontractfigure}}
\end{figure}

The first observation is that, given a specific $P_3$ in a planar graph, we may ``uncontract'' the middle vertex so that the new vertices are adjacent and both adjacent to the end-points of the $P_3$.
We omit a proof since a straight-line drawing makes it clear that the operation is valid, as demonstrated in \Cref{uncontractfigure}.
\begin{obs}[Uncontraction]\label{uncontraction}
    Let $G$ be a planar graph drawn in the plane with straight edges and fix a path on $3$ vertices $xvy$.
    Label the neighbors of $v$ as $x,u_1,\dots,u_k,y,w_1,\ldots,w_{\ell}$ in some cyclic order.
    The operation of \emph{uncontracting along the path $xvy$} creates a new planar graph $G'$ in which the vertex $v$ is replaced by the adjacent vertices $v_1$ and $v_2$, where $N(v_1)=\{x,u_1,\ldots,u_k,y,v_2\}$ and $N(v_2)=\{y,w_1,\ldots,w_\ell,x,v_1\}$.
\end{obs}

Since contracting an edge into a single vertex preserves planarity and we just observed that one can ``uncontract'' a vertex to create a new edge while preserving planarity, we can perform these two operations in sequence: first contracting an edge and then uncontracting the resulting vertex.
See \Cref{uncontractfigure} for a demonstration of this operation, which we dub ``contraction--uncontraction''.

\begin{obs}[Contraction--Uncontraction]\label{contract-uncontract}
    Let $G$ be a planar graph drawn in the plane with straight edges and fix a path on $4$ vertices $xuvy$.
    The operation of \emph{contraction--uncontraction along the path $xuvy$} creates a new graph $G'$ by first contracting the edge $uv$ into a single vertex $(uv)$ and then uncontracting along the $3$-path $x(uv)y$ to recover the vertices $u$ and $v$.
    $G'$ has the same vertex-set as does $G$ and, due to \Cref{uncontraction}, it additionally satisfies:
    \begin{itemize}
        \item $G'$ is planar, and
        \item $uv$ is an edge and both $u$ and $v$ are adjacent to both $x$ and $y$, and
        \item $N_{G'}(u)\cup N_{G'}(v)=N_G(u)\cup N_G(v)$.
    \end{itemize}
\end{obs}

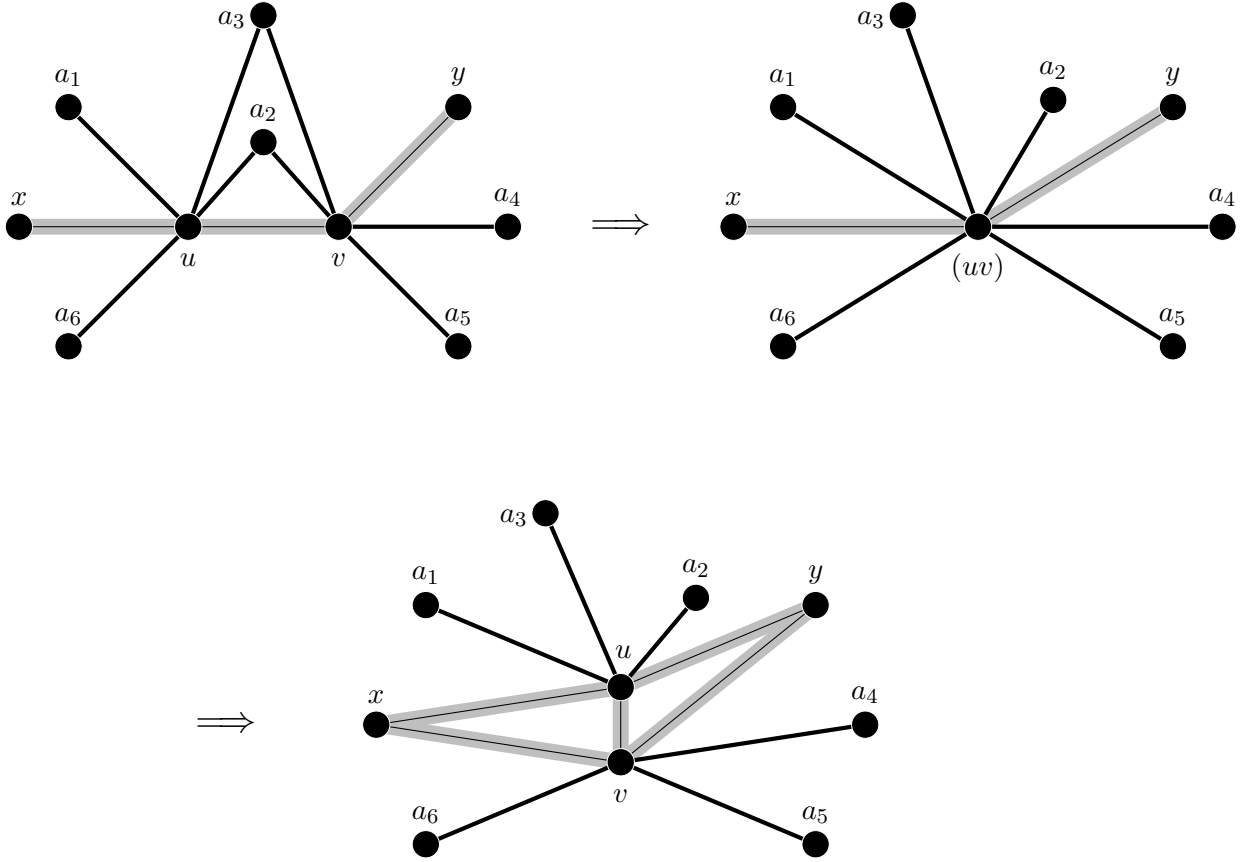
\begin{figure}[ht]
    \begin{center}
        \begin{tikzpicture}
            \def\rad{2.25}; 
            \def\mid{1.0}; 
            \def\spc{1.5}; 
            \coordinate (cen;1) at (-\mid-\rad-\spc,0);
            \coordinate (cen;2) at (\mid+\rad+\spc,0);
            \coordinate (cen;3) at (0,-\mid-2.5*\rad);
            \foreach\i in {1,2,3}{
                \coordinate (x;\i) at ($(cen;\i)-(\mid,0)+(180:\rad)$);
                \coordinate (y;\i) at ($(cen;\i)+(\mid,0)+(45:\rad)$);
                \coordinate (a1;\i) at ($(cen;\i)-(\mid,0)+(135:\rad)$);
                \coordinate (a4;\i) at ($(cen;\i)+(\mid,0)+(0:\rad)$);
                \coordinate (a5;\i) at ($(cen;\i)+(\mid,0)+(-45:\rad)$);
                \coordinate (a6;\i) at ($(cen;\i)-(\mid,0)+(-135:\rad)$);
            }
            \coordinate (a2;1) at ($(cen;1)+(90:0.5*\rad)$);
            \coordinate (a3;1) at ($(cen;1)+(90:1.25*\rad)$);
            \foreach\i in {2,3}{
                \coordinate (a2;\i) at ($(cen;\i)+(\mid,0)+(90:0.75*\rad)$);
                \coordinate (a3;\i) at ($(cen;\i)-(\mid,0)+(90:1.25*\rad)$);
            }
            \coordinate (u;1) at ($(cen;1)-(\mid,0)$);
            \coordinate (v;1) at ($(cen;1)+(\mid,0)$);
            \coordinate (u;2) at (cen;2);
            \coordinate (v;2) at (cen;2);
            \coordinate (u;3) at ($(cen;3)+(0,\mid/2)$);
            \coordinate (v;3) at ($(cen;3)-(0,\mid/2)$);
            \begin{pgfonlayer}{fore}
                \foreach\i in {1,2,3}{
                    \node[vertex,label={[yshift=-1pt]$x$}] at (x;\i) {};
                    \node[vertex,label={[yshift=-1pt]$y$}] at (y;\i) {};
                    \foreach\j in {1,2,4,5,6}{
                        \node[vertex,label={[yshift=-1pt]$a_{\j}$}] at (a\j;\i) {};
                    }
                    \node[vertex,label={[xshift=-12pt,yshift=-14pt]$a_3$}] at (a3;\i) {};
                }
                \node[vertex,label={[yshift=-24pt]$u$}] at (u;1) {};
                \node[vertex,label={[yshift=-24pt]$v$}] at (v;1) {};
                \node[vertex,label={[yshift=-30pt]$(uv)$}] at (u;2) {};
                \node[vertex,label={[xshift=1pt,yshift=2pt]$u$}] at (u;3) {};
                \node[vertex,label={[yshift=-24pt]$v$}] at (v;3) {};
                \node at (0,0) {\Large$\Longrightarrow$};
                \node at ($(cen;3)-(3*\mid+\rad,0)$) {\Large$\Longrightarrow$};
            \end{pgfonlayer}
            \begin{pgfonlayer}{main}
                \foreach\i in {1,2,3}{
                    \foreach\j in {1,2,3}{
                        \draw [edge] (u;\i) -- (a\j;\i);
                    }
                    \foreach\j in {4,5}{
                        \draw [edge] (v;\i) -- (a\j;\i);
                    }
                }
                \draw [edge] (v;1) -- (a2;1);
                \draw [edge] (v;1) -- (a3;1);
                \draw [edge] (u;1) -- (a6;1);
                \foreach\i in {2,3}{
                    \draw [edge] (v;\i) -- (a6;\i);
                }
            \end{pgfonlayer}
            \begin{pgfonlayer}{main}
                \foreach\i in {1,3}{
                    \draw [thedge] (u;\i) -- (v;\i);
                }
                \foreach\i in {1,2,3}{
                    \draw [thedge] (u;\i) -- (x;\i);
                    \draw [thedge] (v;\i) -- (y;\i);
                }
                \draw [thedge] (u;3) -- (y;3);
                \draw [thedge] (v;3) -- (x;3);
            \end{pgfonlayer}
            \begin{pgfonlayer}{back}
                \foreach\i in {1,3}{
                    \draw [spedge] (u;\i) -- (v;\i);
                }
                \foreach\i in {1,2,3}{
                    \draw [spedge] (u;\i) -- (x;\i);
                    \draw [spedge] (v;\i) -- (y;\i);
                }
                \draw [spedge] (u;3) -- (y;3);
                \draw [spedge] (v;3) -- (x;3);
            \end{pgfonlayer}
        \end{tikzpicture}
    \end{center}
\caption{Contraction--uncontraction along the path $xuvy$.\label{contractuncontractfigure}}
\end{figure}

In \Cref{merge-onemors}, we show that, in a tumor graph, if a vertex in $\subs$ has exactly one $\doms$ neighbor (that is, in $\bigsqcup_{x\in\doms}\subs_x$), then under certain conditions we can find a graph with at least as many good cycles with one fewer vertex in $\bigsqcup_{x\in\doms}\subs_x$.
This is the key ingredient necessary to ``clean'' a tumor graph and is accomplished by first contracting an edge and then uncontracting the resulting vertex.

\begin{lemma}[Contraction--Uncontraction]\label{merge-onemors}
    Let $(G;\doms,\subs)$ be a planar tumor graph and fix an edge $uv\in E(G)$ with $u\in\onemor x$ for some $x\in\doms$.
    If
    \begin{enumerate}[label=(\roman*)]
        \item $v\in\onemor y$ for some $y\in\doms\setminus\{x\}$, or \label{merge:onemor}
        \item $v\in\tumor xy$ for some $y\in\doms\setminus\{x\}$ and $N(u)\cap N(v)\subseteq \{x\}\cup \tumor xy$, then \label{merge:tumor}
    \end{enumerate}
    the (planar) graph $G'$ formed by contraction--uncontraction along the path $xuvy$ satisfies
    \[
        \cyc{2m+1}{G';\doms,\subs}\geq\cyc{2m+1}{G;\doms,\subs}.
    \]
\end{lemma}
\begin{proof}
    The result of performing a contraction--uncontraction operation along the path $xuvy$ adds the edges $xv$ and $uy$ (should they not already exist) and perhaps ``scrambles'' the other neighbors of $u$ and $v$.

    Every good cycle will either have one $\subs\subs$ edge and otherwise alternate between vertices in $\doms$ and vertices in $\subs$ or will have one $\doms\doms$ edge and otherwise alternate between vertices in $\doms$ and vertices in $\subs$.

    We will classify the good cycles of both $G$ and $G'$ according to the vertex in the $\subs\subs$ edge (if it exists) that is neither $u$ nor $v$. That is,
    \begin{itemize}
        \item The set $\mathcal{C}_{\emptyset}(G)$ is the set of all good cycles in $G$ that have no $\subs\subs$ edge or for which the $\subs\subs$ edge contains neither $u$ nor $v$ or for which the $\subs\subs$ edge is $uv$ and contains the path $xuvy$. The set $\mathcal{C}_{\emptyset}(G')$ is similarly defined.
        \item The set $\mathcal{C}^{*}(G')$ is the set of all good cycles in $G'$ for which the $\subs\subs$ edge is $uv$ and contains the path $xvuy$. Note that $xvuy$ is not a path in $G$ because $u\in\subs_x$.
        \item For any $w\not\in\{u,v\}$, the set $\mathcal{C}_{w}(G)$ is the set of all good cycles in $G$ for which the $\subs\subs$ edge is either $uw$ or $vw$. The set $\mathcal{C}_{w}(G')$ is similarly defined.
    \end{itemize}

    \noindent\textbf{Case~\labelcref{merge:onemor}.}

    The cycles in $\mathcal{C}_{\emptyset}(G)$ are unchanged after contraction--uncontraction and so map to themselves in $\mathcal{C}_{\emptyset}(G')$. Note that if $uv$ is the $\subs\subs$ edge then the cycle must contain the path $xuvy$.

    For $\mathcal{C}_{w}(G)$, we create a map from $\mathcal{C}_w(G)$ to $\mathcal{C}_w(G')$ according to the subgraph induced by $\{u,v,w\}$.
    \begin{itemize}
        \item If $\{u,v,w\}$ induces a path $vuw$, then the cycle must contain a path $xuwz$ for some $z\in\doms\setminus\{x\}$. Depending on whether $w$ is a neighbor of $u$ or $v$ in $G'$, keep the cycle with $xuwz$ or replace that path with $xvwz$. This is a one-to-one map.
        \item If $\{u,v,w\}$ induces a path $uvw$, then the cycle contains the path $yvwz$ for some $z\in\doms\setminus\{y\}$. We either keep the aforementioned path or replace $v$ with $u$. This is a one-to-one map.
        \item If $\{u,v,w\}$ induces a triangle, there are two types of cycles in this case, those that contain the path $xuw$ and those that contain the path $wvy$. Replace these two paths with either (a) the two paths $xuw,wuy$ or (b) the two paths $xvw,wvy$. This is a two-to-two map.
    \end{itemize}

    \noindent\textbf{Case~\labelcref{merge:tumor}.}

    The cycles in $\mathcal{C}_{\emptyset}(G)$ are unchanged after contraction--uncontraction and so map to themselves in $\mathcal{C}_{\emptyset}(G')$. Note that if $uv$ is the $\subs\subs$ edge then the cycle must contain the path $xuvy$.

    For $\mathcal{C}_{w}(G)$, we make a map from $\mathcal{C}_w(G)$ to $\mathcal{C}_w(G')\cup\mathcal{C}^{*}(G')$ according to the subgraph induced by $\{u,v,w\}$.
    \begin{itemize}
        \item If $\{u,v,w\}$ induces a path $vuw$, then the cycle must contain a path $xuwz$ for some $z\in\doms\setminus\{x\}$. Depending on whether $w$ is a neighbor of $u$ or $v$ in $G'$, keep the cycle with $xuwz$ or replace that path with $xvwz$. This is a one-to-one map.
        \item If $\{u,v,w\}$ induces a path $uvw$, then the cycle either contains the path $xvwz$ for some $z\in\doms\setminus\{x\}$ or contains the path $yvwz$ for some $z\in\doms\setminus\{y\}$. In either case we either keep the aforementioned path or replace $v$ with $u$. This is a one-to-one map.
        \item If $\{u,v,w\}$ induces a triangle, this is the unique $w$ by \Cref{few_in_tumor}. So, there are three types of cycles in this case, those that contain the path $xuwy$, those that contain the path $xwvy$ and those that contain the path $xvwy$. Replace these three paths with either (a) the three paths $xuwy,xwuy,xvuy$ or (b) the three paths $xvwy,xwvy,xvuy$. This is a three-to-three map.\qedhere
    \end{itemize}
\end{proof}

\subsubsection{Cleaning Stage I}
\label{tumor:onemor}

\begin{lemma}[Stage I]\label{stageone}
    Let $(G;\doms,\subs)$ be a planar tumor graph.
    There is another planar tumor graph $(G';\doms,\subs')$ satisfying:
    \begin{enumerate}[label=(\roman*)]
        \item $\subs'\subseteq\subs$, and
        \item $\subs'_{\varnothing}=\varnothing$, and \label{stageone:nomor}
        \item $\bigsqcup_{x\in\doms}\subs'_x$ is an independent set, and \label{stageone:onemor}
        \item $\cyc{2m+1}{G';\doms,\subs'}\geq \cyc{2m+1}{G;\doms,\subs}$.
    \end{enumerate}
\end{lemma}
We say that a planar tumor graph $(G';\doms,\subs')$ with \cref{stageone:nomor,stageone:onemor} is a \emph{Stage I graph}.
\begin{proof}
    We repeatedly modify the graph $G$ until it has the desired properties.

    To begin, set $\subs^{*}=\subs\setminus\nomor$ and let $G^{*}$ be the graph where we remove all vertices within $\nomor$ and remove all edges within $\onemor x$ for each $x\in\doms$.
    Certainly $(G^{*};\doms,\subs^{*})$ is still a planar tumor graph and also the set of good cycles remains unchanged since none of these cycles can use any of the deleted vertices and edges.

    Now, we define $G'$ by repeating the following: while there is an edge $uv$ with $u\in\onemor x$ and $v\in\onemor y$ for some $x\neq y\in\doms$, perform contraction--uncontraction along the path $xuvy$.
    Each time we perform such a contract--uncontract operation, the resulting graph is planar and the size of $\bigsqcup_{x\in\doms}\onemor x$ strictly decreases. So, eventually this process terminates and we have that, in the resulting $G'$, $\bigsqcup_{x\in\doms}\onemor x$ is an independent set.
    Furthermore, setting $\subs'=\subs^{*}$, it is the case that $\cyc{2m+1}{G';\doms,\subs'}\geq\cyc{2m+1}{G^{*};\doms,\subs^{*}}$ by \Cref{merge-onemors}.
\end{proof}

\subsubsection{Cleaning Stage II}
\label{tumor:onemortumor}

First, we observe that removing few edges within $G[\subs]$ results in a negligible reduction in the number of good cycles.

\begin{prop}\label{removeSS}
    Fix $m\geq 2$, let $(G;\doms,\subs)$ be a planar tumor graph and fix any $R\subseteq E(G[\subs])$.
    \[
        \cyc{2m+1}{G;\doms,\subs}\leq\cyc{2m+1}{G-R;\doms,\subs}+O\bigl(\abs R n^{m-1}\bigr).
    \]
\end{prop}
\begin{proof}
    We must count the number of good cycles $(v_1,\dots,v_{2m+1})$ with, say, $v_{2m}v_{2m+1}\in R$.
    Since $v_{2m},v_{2m+1}\in\subs$, we know that $v_{2m-1}\in\doms$.
    Since each member of $\subs$ has at most two neighbors in $\doms$, this yields at most $4\abs R$ many choices for the triple $(v_{2m-1},v_{2m},v_{2m+1})$.
    Then \Cref{pathasymptotics} tells us that there are at most $2\numb(G,P_{2m-2})\leq O(n^{m-1})$ many choices for the path $(v_1,\dots,v_{2m-2})$, which proves the claim.
\end{proof}

We now apply this observation along with contraction--uncontraction to clean a planar tumor graph further.

Recall that tumors $\subs_{xy}$ and $\subs_{zw}$ are distinct if $xy\neq zw$.
\begin{lemma}[Stage II]\label{stagetwo}
    Let $(G;\doms,\subs)$ be a Stage I planar tumor graph on $n$ vertices.
    There is another Stage I planar tumor graph $(G';\doms,\subs)$ that additionally satisfies:
    \begin{enumerate}[label=(\roman*)]
        \item There are no edges between distinct tumors. \label{stagetwo:distinct}
        \item If $v\in\onemor x$ for some $x\in\doms$, then $v$ has at most one neighbor within $\bigsqcup_{yz\in{\doms\choose 2}}\tumor yz$.
            Furthermore, if $v$ has a neighbor within $\tumor yz$, then $x\notin \{y,z\}$. \label{stagetwo:onetwo}
        \item $\cyc{2m+1}{G;\doms,\subs}\leq\cyc{2m+1}{G';\doms,\subs}+O(\abs\doms n^{m-1})$.
    \end{enumerate}
\end{lemma}
We say that a Stage I tumor graph $(G';\doms,\subs)$ with \cref{stagetwo:distinct,stagetwo:onetwo} is a \emph{Stage II graph}.
\begin{proof}
    Define $U\subseteq\bigsqcup_{x\in\doms}\onemor x$ to be the set of vertices $v$ with any of the following properties:
    \begin{itemize}
        \item $v$ has at least three neighbors within $\bigsqcup_{xy\in{\doms\choose 2}}\tumor xy$, or
        \item $v$ has neighbors in distinct tumors $\tumor xy$ and $\tumor zw$, or
        \item $v\in\onemor x$ for some $x\in\doms$ and $v$ has two neighbors within $\tumor yz$ for some $\{y,z\}\not\ni x$.
    \end{itemize}
    Denote by $R$ the set of all edges between $U$ and $\bigsqcup_{xy\in{\doms\choose 2}}\tumor xy$ and set $G_1=G-R$.
    According to \Cref{onemor-tumor}, we have $\abs R\leq O(\abs\doms)$ and so, since $R\subseteq E(G[\subs])$, \Cref{removeSS} implies that
    \[
        \cyc{2m+1}{G;\doms,\subs}\leq\cyc{2m+1}{G_1;\doms,\subs}+O(\abs\doms n^{m-1}).
    \]
    Now, certainly $G_1$ is still planar and $(G_1;\doms,\subs)$ is still a Stage I graph.
    Additionally, if $v\in\onemor x$ for some $x\in\doms$, then
    \begin{itemize}
        \item $v$ has at most two neighbors within $\bigsqcup_{yz\in{\doms\choose 2}}\tumor yz$, and
        \item If $v$ does have two neighbors, then both neighbors reside within the same $\tumor xy$ for some $y\in\doms$.
    \end{itemize}

    Now, we define $G_2$ by repeating the following: while there is an edge $uv$ with $u\in\onemor x$ and $v\in\tumor xy$ for some $x\neq y\in\doms$, perform contraction--uncontraction along the path $xuvy$.
    Each time we perform such a contract--uncontract operation, the size of $\bigsqcup_{x\in\doms}\onemor x$ strictly decreases and so eventually this process terminates, resulting in the planar graph $G_2$.
    Thus, $G_2$ has the property that if $v\in\onemor x$ for some $x\in\doms$, then $v$ has at most one neighbor within $\bigsqcup_{yz\in{\doms\choose 2}}\tumor yz$ and that if $v$ has a neighbor within $\tumor yz$, then $x\notin yz$.
    Now, recalling that $G_1$ was a Stage I graph, at no point in this process do we introduce any new edges incident to $\bigsqcup_{x\in\doms}\onemor x$.
    Thus, if we perform a contract--uncontract operation along the path $xaby$ where $a\in\onemor x$ and $b\in\tumor xy$, then $N(a)\subseteq\{x\}\cup\tumor xy$ at this point.
    In particular, \Cref{merge-onemors} implies that $\cyc{2m+1}{G_1;\doms,\subs}\leq\cyc{2m+1}{G_2;\doms,\subs}$.

    Finally, we form $G'$ by removing all edges between $\tumor xy$ and $\tumor zw$ for all $xy\neq zw\in{\doms\choose 2}$.
    According to \Cref{tumor-tumor}, there are at most $O(\abs\doms)$ many edges of this form and so \Cref{removeSS} implies that $\cyc{2m+1}{G_2;\doms,\subs}\leq\cyc{2m+1}{G';\doms,\subs}+O(\abs\doms n^{m-1})$, which concludes the proof.
\end{proof}

\subsubsection{Cleaning Stage III}
\label{tumor:tumors}

Recall that a tumor graph $(G;\doms,\subs)$ is benign if whenever $uv\in E(G[\subs])$, then $u,v\in\tumor xy$ for some $xy\in{\doms\choose 2}$.
Observe that any tumor graph that is benign is also Stage II.
The only difference between a Stage II and a benign tumor graph is that a Stage II tumor graph can contain edges of the form $uv$ where $u\in\onemor x$ and $v\in\tumor yz$ provided that $x\notin\{y,z\}$.
Thus, the last step needed to prove \Cref{cleaning} is to control all edges of this form.

\begin{prop}\label{newdoms}
    Let $(G;\doms,\subs)$ be a planar tumor graph.
    If $\mathcal Z$ denote the set of all vertices $z\in\subs$ such that $z\in\tumor xy$ and $z$ has some neighbor within $\onemor w$ for some $x,y,w\in\doms$ with $w\notin\{x,y\}$, then $\abs{\mathcal Z}\leq 2\abs\doms$.
\end{prop}
\begin{proof}
    Let $G'$ be the graph formed from $G$ by contracting all edges of the form $wu$ for $w\in\doms,u\in\onemor w$; note that $G'$ is still planar.
    Within the bipartite graph $G'[\mcal Z,\doms]$, each vertex in $\mcal Z$ has at least three neighbors; thus $\abs{\mathcal Z}\leq 2\abs\doms$ by \Cref{bipartite}.
\end{proof}

\begin{lemma}[Stage III]\label{stagethree}
    Let $(G;\doms,\subs)$ be a Stage II planar tumor graph on $n$ vertices.
    There is another planar tumor graph $(G';\doms',\subs')$ such that
    \begin{enumerate}[label=(\roman*)]
        \item $(G';\doms',\subs')$ is benign, and
        \item $\abs{\doms'}\leq 3\abs\doms$, and
        \item $\cyc{2m+1}{G;\doms,\subs}\leq\cyc{2m+1}{G';\doms',\subs'}+O(\abs\doms n^{m-1})$.
    \end{enumerate}
\end{lemma}

Note that applying \Cref{stageone,stagetwo,stagethree} consecutively yields \Cref{cleaning}.

\begin{proof}
    Let $\mcal Z$ be as in \Cref{newdoms}.
    Now, let $R$ denote the set of edges between $\mcal Z$ and $\bigsqcup_{xy\in{\doms\choose 2}}\tumor xy$ and set $G'=G-R$.
    Since $G$ is a Stage II graph, if $z\in\mcal Z\cap\tumor xy$, then the only neighbors of $z$ are in $\{x,y\}\sqcup\tumor xy\sqcup\bigsqcup_{w\in\doms}\onemor w$.
    In particular, any $z\in\mcal Z$ has at most two neighbors within $\bigsqcup_{xy\in{\doms\choose 2}}\tumor xy$ due to \Cref{few_in_tumor}.
    Thus, $\abs R\leq 2\abs{\mcal Z}\leq 4\abs\doms$ by additionally applying \Cref{newdoms}.
    Then, since $R\subseteq E(G[\subs])$, \Cref{removeSS} tells us that
    \[
        \cyc{2m+1}{G;\doms,\subs}\leq\cyc{2m+1}{G';\doms,\subs}+O(\abs\doms n^{m-1}).
    \]

    Observe that if $u\in\subs$, $z\in\mcal Z$ with $uz\in E(G')$, then $u\in\onemor x$ for some $x\in\doms$ and $z\notin\tumor xy$ for any $y\in\doms$.
    Additionally, since $G$ was Stage II, if this is the case then $z$ is the unique neighbor of $u$ within $\mcal Z$.
    Now, set $\doms'=\doms\sqcup\mcal Z$ and $\subs'=\subs\setminus\mcal Z$.
    Therefore, $u\in\tumor xz'$ in $(G',\doms',\subs')$; in particular, $(G';\doms',\subs')$ is a planar tumor graph.
    Consequently, it is quick to observe that this tumor graph is both Stage II and has no edges between $\tumor xy$ and $\tumor zw$ for any distinct $xy, zw\in{\doms'\choose 2}$.

    Thus, $(G';\doms',\subs')$ is benign. Moreover, $\doms'=\doms\sqcup\mcal Z$ and so $\abs{\doms'}\leq \abs\doms+\abs{\mcal Z}\leq 3\abs\doms$.
    \medskip

    Next, consider a good cycle in $(G';\doms,\subs)$ which does not exist in $(G';\doms',\subs')$.
    Any such cycle contains at least one member of $\mcal Z=\subs\cap\doms'$ which appears in a path of the form $\doms\mcal Z\doms$ or of the form $\doms\mcal Z\mcal Z\doms$, otherwise there would be only one member of $\mcal Z$ which would have to be in a path of the form $\doms\mcal Z\subs'\doms$, which is still good in $(G';\doms',\subs')$.

    For cycles with the pattern $\doms\mcal Z\doms$, the number of ways to choose that path is $2\abs{\mcal Z}\leq O(\abs\doms)$ and the number of ways to choose the remaining path is $2\numb(G',P_{2m-2})$.
    By \Cref{pathasymptotics}, $\numb(G',P_{2m-2}) \leq O(n^{m-1})$, and so the total number of such cycles is bounded above by $O(\abs\doms n^{m-1})$.

    For cycles with the pattern $\doms\mcal Z\mcal Z\doms$, the $\mcal Z\mcal Z$ edge must be in the same tumor $\tumor xy$.
    Thus, the number of ways to choose the $\doms\mcal Z\mcal Z\doms$ piece is at most $4e(G[\mcal Z])\leq 12\abs{\mcal Z}\leq O(\doms)$ because $G[\mcal Z]$ is planar. The number of ways to choose the remaining path is $\numb(G',P_{2m-3})\leq O(n^{m-1})$ by \Cref{pathasymptotics}, and so the total number of such cycles is bounded above by $O(\abs\doms n^{m-1})$.

    As a result,
    \[
        \cyc{2m+1}{G';\doms,\subs}\leq\cyc{2m+1}{G';\doms',\subs'}+O(\abs\doms n^{m-1}),
    \]
    which concludes the proof.
\end{proof}

\subsection{Proof of \Cref{reducetomax}: Reduction to maximum likelihood estimators}
\label{sec:reducetomax}

For a set $X$ and a positive integer $k$, we write $(X)_{k}$ to indicate the set of all tuples $(x_1,\dots,x_k)\in X^k$ with $x_1,\dots,x_k$ distinct.
This notation mirrors that of the falling-factorial.
\medskip

We begin by constructing an edge probability measure $\mu$ on the clique with vertex set $\doms$ where
\[
    \mu(xy)\eqdef\frac{\abs {\tumor xy}}{\sum_{zw\in\binom{\doms}{2}}\abs {\tumor zw}}, \qquad\qquad\text{for all }xy\in{\doms\choose 2}.
\]
Since the tumors are disjoint, we know that $\abs{\tumor xy}\leq\mu(xy)\cdot n$.
Furthermore, by \Cref{few_in_tumor}, $e\bigl(G_2[\tumor xy]\bigr)\leq\abs{\tumor xy}\leq \mu(xy)\cdot n$.
\medskip

Next, recall that the good cycles in $(G_2;\doms,\subs)$ alternate between $\doms$ vertices and $\subs$ vertices except for one consecutive pair which can either be $\doms\doms$ or $\subs\subs$.
Let $\mcal C_\subs$ denote those good copies of $C_{2m+1}$ containing an $\subs\subs$ edge and let $\mcal C_\doms$ denote those good copies of $C_{2m+1}$ containing a $\doms\doms$ edge.
Of course, $\cyc{2m+1}{G_2;\doms,\subs}=\abs{\mcal C_\subs}+\abs{\mcal C_\doms}$.
\medskip

Fix a cycle in $\mcal C_\subs$ and label its vertices cyclically as $(v_1,\dots,v_{2m+1})$ so that $v_{2m},v_{2m+1}\in\subs$.
Then $(x_1,x_2,\dots,x_m)=(v_1,v_3,\dots,v_{2m-1})$ has the property that $(x_1,\dots,x_m)\in(\doms)_m$ and $v_{2i}\in\tumor{x_i}{x_{i+1}}$ for all $i\in[m-1]$.
Furthermore, since $(G_2;\doms,\subs)$ is refined, we know that $v_{2m}v_{2m+1}\in E\bigl(G_2[\tumor{x_m}{x_1}]\bigr)$.
Thus, the number of cycles in $\mcal C_\subs$ which yield the tuple $(x_1,\dots,x_m)\in(\doms)_m$ is precisely
\[
    \biggl(\prod_{i=1}^{m-1} \abs{\tumor{x_i}{x_{i+1}}}\biggr)\cdot 2 e\bigl(G_2[\tumor{x_1}{x_m}]\bigr).
\]
Of course, there are two cyclic orderings of the vertices of each of these cycles and so
\begin{align*}
    \abs{\mcal C_\subs} &={1\over 2}\sum_{(x_1,\dots,x_m)\in(\doms)_m}\biggl(\prod_{i=1}^{m-1} \abs{\tumor{x_i}{x_{i+1}}}\biggr)\cdot 2 e\bigl(G_2[\tumor{x_1}{x_m}]\bigr)\\
                        &=\sum_{(x_1,\dots,x_m)\in(\doms)_m}\biggl(\prod_{i=1}^{m-1} \abs{\tumor{x_i}{x_{i+1}}}\biggr)\cdot e\bigl(G_2[\tumor{x_1}{x_m}]\bigr).
\end{align*}
Bringing in the edge probability measure $\mu$, we further bound
\begin{equation}\label{eqn:subsub}
    \abs{\mcal C_\subs} \leq \sum_{(x_1,\dots,x_m)\in(\doms)_m}\biggl(\prod_{i=1}^{m-1}\mu(x_ix_{i+1})\biggr)\cdot\mu(x_mx_1)\cdot n^m
    = \begin{cases}
        \displaystyle 2\sum_{e\in\supp\mu}\mu(e)^2\cdot n^2 & \text{, if }m=2;\\[20pt]
        \displaystyle 2m\cdot\beta(\mu;C_m)\cdot n^m &\text{, if }m\geq 3.
    \end{cases}
\end{equation}
\medskip

Next, fix a cycle in $\mcal C_\doms$ and label its vertices cyclically as $(v_1,\dots,v_{2m+1})$ so that $v_{1},v_{2m+1}\in\doms$.
Then $(x_1,x_2,\dots,x_{m+1})=(v_1,v_3,\dots,v_{2m+1})$ has the property that $(x_1,\dots,x_{m+1})\in(\doms)_{m+1}$ and $v_{2i}\in\tumor{x_i}{x_{i+1}}$ for all $i\in[m]$.
Thus, the number of cycles in $\mcal C_\doms$ which yield the tuple $(x_1,\dots,x_{m+1})\in(\doms)_{m+1}$ is precisely
\[
    \prod_{i=1}^m\abs{\tumor{x_i}{x_{i+1}}},\qquad\text{provided that }x_1x_m\in E(G_2).
\]
Again, there are two cyclic orderings of the vertices of each of these cycles and so
\[
    \abs{\mcal C_\doms} = {1\over 2}\sum_{\substack{(x_1,\dots,x_{m+1})\in(\doms)_{m+1},\\ x_1x_{m+1}\in E(G_2)}} \prod_{i=1}^m\abs{\tumor{x_i}{x_{i+1}}}.
\]
By dropping the requirement that $x_1x_m\in E(G_2)$ and also bringing in the edge probability measure $\mu$, we bound
\begin{align}
    \abs{\mcal C_\doms} &\leq{1\over 2}\sum_{(x_1,\dots,x_{m+1})\in(\doms)_{m+1}}\prod_{i=1}^m\abs{\tumor{x_i}{x_{i+1}}}\nonumber \\
                        &\leq{1\over 2}\sum_{(x_1,\dots,x_{m+1})\in(\doms)_{m+1}}\prod_{i=1}^m\mu(x_ix_{i+1})\cdot n^m\nonumber \\
                        &=\beta(\mu;P_{m+1})\cdot n^m\label{eqn:domdom}
\end{align}

Combining \cref{eqn:subsub,eqn:domdom} finally yields
\begin{align*}
    \cyc{5}{G_2;\doms,\subs} &\leq \biggl(2\sum_{e\in\supp\mu}\mu(e)^2+\beta(\mu;P_3)\biggr)n^2, && \text{and}\\
    \cyc{2m+1}{G_2;\doms,\subs} &\leq \bigl(2m\cdot\beta(\mu;C_m)+\beta(\mu;P_{m+1})\bigr)n^m, && \text{for all $m\geq 3$},
\end{align*}
which concludes the proof of \Cref{reduction}. \hfill\qedsymbol

\section{Proof of \Cref{maxlikelihood}}
\label{maxlikelihood-proof}

We first address the case where $m=2$.
\begin{prop}\label{C5-max}
    If $\mu$ is an edge probability measure, then
    \[
        2\sum_{e\in\supp\mu}\mu(e)^2+\beta(\mu;P_3)\leq 2,
    \]
    with equality if and only if $\abs{\supp\mu}=1$.
\end{prop}
\begin{proof}
    We use the definition of $\beta(\mu;P_3)$:
    \begin{align*}
        2\sum_{e\in\supp\mu}\mu(e)^2+\beta(\mu;P_3)
        &= 2\sum_{e\in\supp\mu}\mu(e)^2+\sum_{ef\in{\supp\mu\choose 2}:\ \abs{e\cap f}=1}\mu(e)\mu(f)\\
        &\leq 2\sum_{e\in\supp\mu}\mu(e)^2+{1\over 2}\sum_{e\in\supp\mu}\sum_{f\in\supp\mu\setminus\{e\}}\mu(e)\mu(f)\\
        &\leq 2\biggl(\sum_{e\in\supp\mu}\mu(e)^2 + \sum_{e\in\supp\mu}\sum_{f\in\supp\mu\setminus\{e\}}\mu(e)\mu(f)\biggr)\\
        &= 2\biggl(\sum_{e\in\supp\mu}\mu(e)\biggr)^2 = 2.
    \end{align*}
    Note that the last inequality is an equality if and only if $\abs{\supp\mu}=1$.
\end{proof}

The remainder of this section is dedicated to bounding $2m\cdot\beta(\mu;C_m)+\beta(\mu;P_{m+1})$ for $m\geq 3$.
\medskip

Many of our arguments focus on the mass of a vertex in an edge probability measure:

\begin{defn}
    Fix an edge probability measure $\mu\in\Delta^K$ for some clique $K$.
    The function $\bar\mu\colon V(K)\to\R$ is defined by
    \[
        \bar\mu(x)\eqdef\sum_{y\in V(K)\setminus\{x\}}\mu(xy).
    \]
\end{defn}

That is, $\bar\mu(x)$ is the probability that an edge sampled from $\mu$ is incident to the vertex $x$, and can be understood as the weighted degree of $x$.
Note that $\sum_{x\in V(K)}\bar\mu(x)=2$ (handshaking lemma).
\medskip

The next lemma, \Cref{regularity}, is a very general statement in the setting of the maximum likelihood graph problems, that establishes regularity conditions for local optimizers. It is a direct consequence of the Karush--Kuhn--Tucker (KKT) conditions (see {\cite[Corollaries 9.6 and 9.10]{guler_opt}}). We will apply the lemma in the case of $k=2$, $H_1=C_m$, and $H_2=P_{m+1}$.

\begin{lemma}\label{regularity}
    Fix a positive integer $m$ and any graphs $H_1,\dots,H_k$, each with $m$ edges.
    Additionally, fix any clique $K$ on at least two vertices.
    For constants $\gamma_1,\dots,\gamma_k$, set
    \[
        \mcal O= \max_{\nu\in\Delta^K}\ \sum_{i=1}^k\gamma_i\cdot\beta(\nu;H_i).
    \]
    If $\mu\in\Delta^K$ achieves $\mcal O$, then
    \begin{align*}
        m\cdot\mcal O\cdot\mu(e) &= \sum_{i=1}^k\gamma_i\sum_{\substack{H\in\cp(K,H_i),\\ E(H)\ni e}}\mu(H), &&\text{for each }e\in E(K); \\
        m\cdot\mcal O\cdot\bar\mu(x) &= \sum_{i=1}^k\gamma_i\sum_{\substack{H\in\cp(K,H_i),\\ V(H)\ni x}}\deg_H(x)\mu(H), &&\text{for each }x\in V(K).
    \end{align*}
\end{lemma}
We quickly remark that the above maximum is indeed achieved since $\Delta^K$ is compact for any clique $K$ and $\beta(\cdot;H)$ is a continuous function on $\Delta^K$.
\begin{proof}
    By definition, we can write
    \[
        \begin{array}{cccl}
            \mcal O
            & = & \max & \displaystyle\sum_{i=1}^k \gamma_i \sum_{H\in\cp(K,H_i)}\biggl(\prod_{e\in E(H)}\nu(e)\biggr)\\[2em]
            & & \text{s.t.} & \displaystyle\sum_{e\in E(K)}\nu(e) = 1\\[2em]
            & & & \displaystyle \nu(e)\geq 0, \qquad\text{for all }e\in E(K).
        \end{array}
    \]
    In particular, we may apply the KKT conditions to this optimization problem to find that if $\mu\in\Delta^K$ achieves $\mcal O$, then there is some fixed $\lambda\in\R$ such that $D(e)=\lambda$ for all $e\in\supp\mu$, where
    \[
        D(e)=\sum_{i=1}^k\gamma_i\sum_{\substack{H\in\cp(K,H_i),\\ E(H)\ni e}}\biggl(\prod_{f\in E(H)\setminus\{e\}}\mu(f)\biggr).
    \]
    Of course, whether or not $e\in\supp\mu$, we always have
    \begin{equation}\label{KKTedge}
        \lambda\cdot\mu(e)=D(e)\cdot\mu(e)=\sum_{i=1}^k\gamma_i\sum_{\substack{H\in\cp(K,H_i),\\E(H)\ni e}}\mu(H).
    \end{equation}
    By then summing over all $e\in E(K)$, we find
    \begin{align*}
        \lambda
        = \sum_{e\in E(K)}\lambda\cdot\mu(e) &= \sum_{e\in E(K)}\sum_{i=1}^k\gamma_i\sum_{\substack{H\in\cp(K,H_i),\\E(H)\ni e}}\mu(H)\\
        &= \sum_{i=1}^k\gamma_i \sum_{H\in\cp(K,H_i)}\mu(H)\sum_{e\in E(K)}\mathbf 1[e\in E(H)]\\
        &= \sum_{i=1}^k\gamma_i \sum_{H\in\cp(K,H_i)}\mu(H)\cdot m = m\cdot\mcal O,
    \end{align*}
    where the penultimate equality follows from the assumption that each $H_i$ has exactly $m$ edges.
    Substituting this value of $\lambda=m\cdot\mcal O$ into \cref{KKTedge} yields the first part of the lemma.

    For the second part of the lemma, we use the first part to find that for any fixed $x\in V(K)$,
    \begin{align*}
        m\cdot\mcal O\cdot\bar\mu(x)
        &= \sum_{y\in V(K)\setminus\{x\}}m\cdot\mcal O\cdot\mu(xy) = \sum_{y\in V(K)\setminus\{x\}}\sum_{i=1}^k\gamma_i\sum_{\substack{H\in\cp(K,H_i),\\ E(H)\ni xy}}\mu(H)\\
        &= \sum_{i=1}^k\gamma_i\sum_{H\in\cp(K,H_i)}\mu(H)\sum_{y\in V(K)\setminus\{x\}}\mathbf 1[xy\in E(H)]\\
        &= \sum_{i=1}^k\gamma_i\sum_{H\in\cp(K,H_i)}\mu(H)\deg_H(x).\qedhere
    \end{align*}
\end{proof}

We primarily use \Cref{regularity} to understand how $\beta(\mu;H)$ changes upon deleting a vertex from $\supp\bar\mu$, which is key to the proof of \Cref{vertex-bound}.
Before we can established \Cref{vertex-bound}, we need a brief, general fact about paths.

\begin{prop}\label{rootedpath}
    Fix a clique $K$ and a vertex $x\in V(K)$.
    For any $\mu\in\Delta^K$ and any integer $m\geq 2$,
    \[
        \sum_{\substack{P\in\cp(K,P_{m+1}),\\ \deg_P(x)=1}}\mu(P)\leq \bar\mu(x)\biggl({1-\bar\mu(x)\over m-1}\biggr)^{m-1}.
    \]
\end{prop}
\begin{proof}
    We prove this by induction on $m$, starting with $m=2$.
    In this case, we have
    \begin{align*}
        \sum_{\substack{P\in\cp(K,P_{m+1}),\\ \deg_P(x)=1}}\mu(P)
        &=\sum_{y\in V(K)\setminus\{x\}}\mu(xy)\sum_{z\in V(K)\setminus\{y,x\}}\mu(yz)=\sum_{y\in V(K)\setminus\{x\}}\mu(xy)\bigl(\bar\mu(y)-\mu(xy)\bigr)\\
        &\leq \sum_{y\in V(K)\setminus\{x\}}\mu(xy)\bigl(1-\bar\mu(x)\bigr)=\bar\mu(x)\bigl(1-\bar\mu(x)\bigr),
    \end{align*}
    where the inequality follows from the fact that $\bar\mu(x)+\bar\mu(y)\leq 1+\mu(xy)$.

    Now suppose that $m\geq 3$.
    Observe that if $\bar\mu(x)=1$, then the inequality trivially holds since there are no positive-mass copies of $P_m$ emanating from $x$.
    Therefore, we may suppose that $\bar\mu(x)<1$ and define a new probability mass $\nu\in\Delta^K$ by effectively deleting the edges incident to $x$:
    \[
        \nu(e)\eqdef{1\over 1-\bar\mu(x)}\cdot\begin{cases}
            0 & \text{if }e\ni x\\
            \mu(e) & \text{otherwise}.
        \end{cases}
    \]
    Then we have
    \begin{align*}
        \sum_{\substack{P\in\cp(K,P_{m+1}),\\ \deg_P(x)=1}}\mu(P)
        &= \sum_{y\in V(K)\setminus\{x\}}\mu(xy)\sum_{\substack{P\in\cp(K,P_m),\\ \deg_P(y)=1,\ P\not\ni x}}\mu(P)\\
        &=\sum_{y\in V(K)\setminus\{x\}}\mu(xy)\cdot\bigl(1-\bar\mu(x)\bigr)^{m-1}\sum_{\substack{P\in\cp(K,P_m),\\ \deg_P(y)=1}}\nu(P)\\
        &\leq\sum_{y\in V(K)\setminus\{x\}}\mu(xy)\cdot\bigl(1-\bar\mu(x)\bigr)^{m-1}\cdot\bar\nu(y)\biggl({1-\bar\nu(y)\over m-2}\biggr)^{m-2}.
    \end{align*}
    Applying the inequality of arithmetic and geometric means (AM--GM inequality) to the pieces of the above expression involving $\bar\nu(y)$ then yields
    \begin{align*}
        \sum_{\substack{P\in\cp(K,P_{m+1}),\\ \deg_P(x)=1}}\mu(P)
        &\leq\sum_{y\in V(K)\setminus\{x\}}\mu(xy)\cdot\bigl(1-\bar\mu(x)\bigr)^{m-1}\cdot\biggl({\bar\nu(y)+(m-2){1-\bar\nu(y)\over m-2}\over m-1}\biggr)^{m-1}\\
        &=\sum_{y\in V(K)\setminus\{x\}}\mu(xy)\cdot\biggl({1-\bar\mu(x)\over m-1}\biggr)^{m-1}=\bar\mu(x)\biggl({1-\bar\mu(x)\over m-1}\biggr)^{m-1}.\qedhere
    \end{align*}
\end{proof}

We now use the above facts to derive an inequality on the vertex-masses in an optimal measure.

\begin{lemma}\label{vertex-bound}
    Fix an integer $m\geq 3$ and fix a clique $K$ on at least $m$ vertices.
    Set
    \[
        \mcal O=\max_{\nu\in\Delta^K}\bigl(2m\cdot\beta(\nu,C_m)+\beta(\nu;P_{m+1})\bigr).
    \]
    If $\mu\in\Delta^K$ achieves $\mcal O$, then
    \[
        {m\over 2}\bar\mu(x)+\bigl(1-\bar\mu(x)\bigr)^m+{1\over 2\mcal O}\cdot\bar\mu(x)\biggl({1-\bar\mu(x)\over m-1}\biggr)^{m-1}\geq 1,\qquad\text{for all } x\in V(K).
    \]
\end{lemma}
\begin{proof}
    Note that $\mcal O>0$ since the uniform distribution on $K$ contains positive-mass copies of $C_m$ since $K$ has at least $m$ vertices.

    Fix any $x\in V(K)$ and note that $\bar\mu(x)<1$.
    We define a new probability mass $\nu\in\Delta^K$ by effectively deleting the edges incident to $x$:
    \[
        \nu(e) \eqdef {1\over 1-\bar\mu(x)}\cdot\begin{cases}
            0 & \text{if }e\ni x\\
            \mu(e) & \text{otherwise}.
        \end{cases}
    \]
    Since $\nu\in\Delta^K$ and $\mcal O$ is the optimal value, we bound
    \begin{align*}
        \mcal O
        &\geq 2m\cdot\beta(\nu;C_m)+\beta(\nu;P_{m+1})\\
        &=2m\sum_{C\in\cp(K,C_m)}\nu(C)+\sum_{P\in\cp(K,P_{m+1})}\nu(P)\\
        &=2m\sum_{\substack{C\in\cp(K,C_m),\\ V(C)\not\ni x}}{\mu(C)\over\bigl(1-\bar\mu(x)\bigr)^m} + \sum_{\substack{P\in\cp(K,P_{m+1}),\\ V(P)\not\ni x}}{\mu(P)\over\bigl(1-\bar\mu(x)\bigr)^m}\\
        &= {1\over\bigl(1-\bar\mu(x)\bigr)^m}\cdot\biggl(\mcal O-2m\sum_{\substack{C\in\cp(K,C_{m}),\\ V(C)\ni x}}\mu(C)-\sum_{\substack{P\in\cp(K,P_{m+1}),\\ V(P)\ni x}}\mu(P)\biggr).
    \end{align*}
    Rearranging this expression yields
    \begin{equation}\label{eqn:CandP}
        1 \leq (1-\bar\mu(x))^m +{1\over\mcal O}\biggl(2m\sum_{\substack{C\in\cp(K,C_m),\\ V(C)\ni x}}\mu(C)+\sum_{\substack{P\in\cp(K,P_{m+1}), \\V(P)\ni x}}\mu(P)\biggr).
    \end{equation}

    Next, \Cref{regularity} tells us that
    \begin{align*}
        m\cdot\mcal O\cdot\bar\mu(x)
        &=2m\sum_{\substack{C\in\cp(K,C_m),\\ V(C)\ni x}}\mu(C)\deg_C(x)+\sum_{\substack{P\in\cp(K,P_{m+1}),\\ V(P)\ni x}}\mu(P)\deg_P(x)\\
        &= 4m\sum_{\substack{C\in\cp(K,C_m),\\ V(C)\ni x}}\mu(C) + \biggl(2\sum_{\substack{P\in\cp(K,P_{m+1}),\\ V(P)\ni x}}\mu(P)-\sum_{\substack{P\in\cp(K,P_{m+1}),\\ \deg_P(x)=1}}\mu(P)\biggr),
    \end{align*}
    and so
    \[
        2m\sum_{\substack{C\in\cp(K,C_m),\\ V(C)\ni x}}\mu(C)+\sum_{\substack{P\in\cp(K,P_{m+1}),\\ V(P)\ni x}}\mu(P)={m\over 2}\cdot\mcal O\cdot\bar\mu(x)+{1\over 2}\sum_{\substack{P\in\cp(K,P_{m+1}),\\ \deg_P(x)=1}}\mu(P).
    \]
    Substituting this expression into \cref{eqn:CandP} and then applying \Cref{rootedpath} finally yields the claim:
    \begin{align*}
        1 &\leq \bigl(1-\bar\mu(x)\bigr)^m+ {m\over 2}\bar\mu(x) +{1\over 2\mcal O}\sum_{\substack{P\in\cp(K,P_{m+1}),\\ \deg_P(x)=1}}\mu(P) \\
          &\leq \bigl(1-\bar\mu(x)\bigr)^m+{m\over 2}\bar\mu(x)+{1\over 2\mcal O}\cdot\bar\mu(x)\biggl({1-\bar\mu(x)\over m-1}\biggr)^{m-1}.\qedhere
    \end{align*}
\end{proof}

We now solve the maximum likelihood question in the case where $m\in\{3,4\}$.

\begin{theorem}\label{exact34}
    For $m\in\{3,4\}$ and any edge probability measure $\mu$,
    \[
        2m\cdot\beta(\mu;C_m)+\beta(\mu;P_{m+1})\leq {2\over m^{m-1}}
    \]
    with equality if and only if $\mu$ is the uniform distribution on $E(C_m)$.
\end{theorem}
\begin{proof}
    Fix a clique $K$ on at least $m$ vertices and set
    \[
        \mcal O = \max_{\nu\in\Delta^K}\bigl(2m\cdot\beta(\nu;C_m)+\beta(\nu;P_{m+1})\bigr).
    \]
    Note that $\mcal O\geq 2/m^{m-1}$ since this is the value achieved by the uniform distribution on a copy of $C_m$, which is a member $\Delta^K$.

    Fix any mass $\mu\in\Delta^K$ which achieves $\mcal O$.
    By \Cref{vertex-bound} and the fact that $\mcal O\geq 2/m^{m-1}$, for each $x\in V(K)$,
    \begin{align}
        1&\leq {m\over 2}\bar\mu(x)+(1-\bar\mu(x))^m+{\bar\mu(x)\over 2\mcal O}\biggl({1-\bar\mu(x)\over m-1}\biggr)^{m-1} \nonumber \\
         &\leq {m\over 2}\bar\mu(x)+(1-\bar\mu(x))^m+{\bar\mu(x)\cdot m^{m-1}\over 4}\biggl({1-\bar\mu(x)\over m-1}\biggr)^{m-1}. \label{eqn:vertbound}
    \end{align}

    Lv et al.~\cite{lv_cycles} proved that $\beta(\mu;C_m)\leq 1/m^m$ with equality if and only if $\mu$ is the uniform distribution on $E(C_m)$.
    Thus, in order to conclude the proof, it suffices to prove that $\abs{\supp\bar\mu}\leq m$ since then $\beta(\mu;P_{m+1})=0$ trivially.

    \noindent\textbf{Case~$m=3$.}
    In this case, \cref{eqn:vertbound} implies that either $\bar\mu(x)=0$ or $\bar\mu(x)>0.57$.
    Thus,
    \[
        2=\sum_{x\in\supp\bar\mu}\bar\mu(x)>\abs{\supp\bar\mu}\cdot 0.57\quad\implies\quad \abs{\supp\bar\mu} < 3.51,
    \]
    and so $\abs{\supp\bar\mu}\leq 3$.
    \medskip

    \noindent\textbf{Case~$m=4$.}
    In this case, \cref{eqn:vertbound} implies that either $\bar\mu(x)=0$ or $\bar\mu(x)>0.41$; therefore, $\abs{\supp\bar\mu}<4.88\implies \abs{\supp\bar\mu}\leq 4$.
\end{proof}

We conclude by establishing a bound on $2m\cdot\beta(\mu;C_m)+\beta(\mu;P_{m+1})$ for all $m\geq 5$ which is tight up to a constant which is independent of $m$.

\begin{theorem}\label{somebound}
    For all $m\geq 5$ and any edge probability measure $\mu$,
    \[
        2m\cdot\beta(\mu;C_m)+\beta(\mu;P_{m+1})\leq {2.6947\over m^{m-1}}
    \]
\end{theorem}
\begin{proof}
    As above, fix a clique $K$ on at least $m$ vertices and set
    \[
        \mcal O = \max_{\nu\in\Delta^K}\bigl(2m\cdot\beta(\mu;C_m)+\beta(\mu;P_{m+1})\bigr).
    \]
    Again, note that $\mcal O\geq 2/m^{m-1}$.

    Fix a mass $\mu\in\Delta^K$ which achieves $\mcal O$ and set
    \[
        s = m\cdot\min_{x\in\supp\bar\mu} \bar\mu(x).
    \]
    By \Cref{vertex-bound} and the fact that $\mcal O\geq 2/m^{m-1}$, we know that
    \begin{align*}
        1 &\leq {m\over 2}{s\over m}+\biggl(1-{s\over m}\biggr)^m+{m^{m-1}\over 4}{s\over m}\biggl({1-{s\over m}\over m-1}\biggr)^{m-1}\\
          &= {s\over 2}+\biggl(1-{s\over m}\biggr)^m + {s\over 4m}\biggl(1-{s\over m}\biggr)^{m-1}\biggl(1+{1\over m-1}\biggr)^{m-1}\\
          &\leq {s\over 2}+e^{-s}+{s\over 4m}e^{1-s+s/m}\leq {s\over 2}+e^{-s}+{s\over 20}e^{1-4s/5},
    \end{align*}
    where the final inequality follows from the fact that $m\geq 5$.
    This implies that $s\geq 1.3644$.
    \medskip

    Next, we bound $\beta(\mu;P_{m+1})$.
    Here, we use \Cref{rootedpath} along with the lower bound of $s\geq 1.3644$ to bound
    \begin{align*}
        \beta(\mu;P_{m+1}) &=\sum_{P\in\cp(K,P_{m+1})}\mu(P) = {1\over 2}\sum_{x\in\supp\bar\mu}\ \sum_{\substack{P\in\cp(K,P_{m+1}),\\ \deg_P(x)=1}}\mu(P) \leq {1\over 2}\sum_{x\in\supp\bar\mu} \bar\mu(x)\biggl({1-\bar\mu(x)\over m-1}\biggr)^{m-1}\\
                           &\leq {1\over 2}\sum_{x\in\supp\bar\mu} \bar\mu(x)\biggl({1-{s\over m}\over m-1}\biggr)^{m-1} =\biggl({1-{s\over m}\over m-1}\biggr)^{m-1} \leq {e^{1-s}\over m^{m-1}}\leq {e^{-0.3644}\over m^{m-1}}.
    \end{align*}

    Finally, we use the fact that $\beta(\mu;C_m)\leq 1/m^m$ to bound
    \[
        2m\cdot\beta(\mu;C_m)+\beta(\mu;P_{m+1}) \leq {2\over m^{m-1}}+{e^{-0.3644}\over m^{m-1}}\leq {2.6947\over m^{m-1}}.\qedhere
    \]
\end{proof}

This concludes the proof of \Cref{maxlikelihood}.

\section{Concluding remarks}

The main question left open by this paper is that of bounding $2m\cdot\beta(\mu;C_m)+\beta(\mu;P_{m+1})$.
\begin{conj}
    For every $m\geq 3$ and every edge probability mass $\mu$,
    \[
        2m\cdot\beta(\mu;C_m)+\beta(\mu;P_{m+1})\leq{2\over m^{m-1}},
    \]
    with equality if and only if $\mu$ is the uniform distribution on $E(C_m)$.
\end{conj}
If true, then
\[
    \numb_\plan(n,C_{2m+1})=2m\biggl({n\over m}\biggr)^m+O(n^{m-1/5})\qquad\text{for all }m\geq 3.
\]

In fact, we believe the stronger bound $2m\cdot\beta(\mu;C_m)+2\cdot\beta(\mu;P_{m+1})\leq 2/m^{m-1}$ to hold, though this has no bearing on the corresponding question of bounding $\numb_\plan(n,C_{2m+1})$.

\paragraph{Even paths.}

In \cite{cox_paths}, Cox and Martin additionally proved a reduction lemma for paths on an odd number of vertices which used many of the same ideas as their reduction lemma for even cycles.
It is natural to wonder if the ideas introduced in this paper can be applied to produce an analogous reduction lemma for even paths.
This is especially motivated by the fact that the conjectured (asymptotic) extremal structure for $\numb_\plan(n,P_{2m})$ is identical to that for $\numb_\plan(n,C_{2m+1})$, namely a balanced blow-up of $C_m$ (see~\cite{ghosh_planarp5}).

\Cref{tumorsamecycles,fewbadcycles} have direct analogues when trying to bound $\numb_\plan(n,P_{2m})$; in fact, the proof that most copies of $P_{2m}$ in a planar tumor graph are ``good'' (contain at most one instance of $\doms\doms$ or $\subs\subs$) is arguably simpler than the proof of \Cref{fewbadcycles}.
Furthermore, there is a direct analogue to \Cref{reducetomax} relating the number of good copies of $P_{2m}$ in a benign planar tumor graph to a maximum likelihood problem, although this maximum likelihood problem is significantly more complex.
Unfortunately, there are major obstructions to proving an analogue of \Cref{cleaning}, the cleaning lemma.
The main operation used in the cleaning lemma is contraction--uncontraction (\Cref{contract-uncontract}) in order to rearrange misbehaving edges.
Our argument that contraction--uncontraction does not decrease the total number of good cycles relied on ``locally rerouting'' the good cycles (\Cref{merge-onemors}).
That is to say, we made no global considerations about the total number of good cycles nor their overall structure.

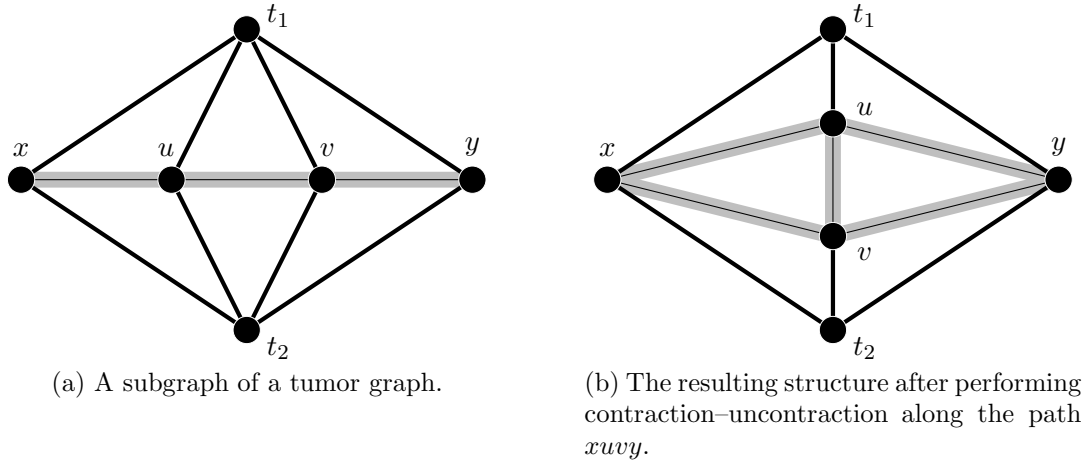
\begin{figure}[ht]
    \centering
    \begin{subfigure}[b]{0.4\textwidth}
        \begin{tikzpicture}
            \useasboundingbox (-3.3,-2.3) rectangle (3.3,2.3);
            \coordinate (x) at (-3,0);
            \coordinate (y) at (3,0);
            \coordinate (u) at (-1,0);
            \coordinate (v) at (1,0);
            \coordinate (t1) at (0,2);
            \coordinate (t2) at (0,-2);
            \begin{pgfonlayer}{fore}
                \node[vertex,label={[yshift=0pt]$x$}] at (x) {};
                \node[vertex,label={[yshift=0pt]$y$}] at (y) {};
                \node[vertex,label={[xshift=-2pt,yshift=0pt]$u$}] at (u) {};
                \node[vertex,label={[xshift=2pt,yshift=0pt]$v$}] at (v) {};
                \node[vertex,label={[xshift=12pt,yshift=-8pt]$t_1$}] at (t1) {};
                \node[vertex,label={[xshift=12pt,yshift=-20pt]$t_2$}] at (t2) {};
            \end{pgfonlayer}
            \begin{pgfonlayer}{main}
                \draw[thedge] (x)--(u)--(v)--(y);
                \draw[edge] (t1)--(x)--(t2);
                \draw[edge] (t1)--(u)--(t2);
                \draw[edge] (t1)--(v)--(t2);
                \draw[edge] (t1)--(y)--(t2);
            \end{pgfonlayer}
            \begin{pgfonlayer}{back}
                \draw[spedge] (x)--(u)--(v)--(y);
            \end{pgfonlayer}
        \end{tikzpicture}
        \caption{\label{fig:path-obstruction:before}A subgraph of a tumor graph.~\\~\\~}
    \end{subfigure}
    \hfil
    \begin{subfigure}[b]{0.4\textwidth}
        \begin{tikzpicture}
            \useasboundingbox (-3.3,-2.3) rectangle (3.3,2.3);
            \coordinate (x) at (-3,0);
            \coordinate (y) at (3,0);
            \coordinate (u) at (0,0.75);
            \coordinate (v) at (0,-0.75);
            \coordinate (t1) at (0,2);
            \coordinate (t2) at (0,-2);
            \begin{pgfonlayer}{fore}
                \node[vertex,label={[yshift=0pt]$x$}] at (x) {};
                \node[vertex,label={[yshift=0pt]$y$}] at (y) {};
                \node[vertex,label={[xshift=12pt,yshift=-6pt]$u$}] at (u) {};
                \node[vertex,label={[xshift=12pt,yshift=-18pt]$v$}] at (v) {};
                \node[vertex,label={[xshift=12pt,yshift=-8pt]$t_1$}] at (t1) {};
                \node[vertex,label={[xshift=12pt,yshift=-20pt]$t_2$}] at (t2) {};
            \end{pgfonlayer}
            \begin{pgfonlayer}{main}
                \draw[edge] (t1)--(u);
                \draw[edge] (v)--(t2);
                \draw[edge] (x)--(t1)--(y);
                \draw[edge] (x)--(t2)--(y);
                \draw[thedge] (u)--(v);
                \draw[thedge] (x)--(u)--(y);
                \draw[thedge] (x)--(v)--(y);
            \end{pgfonlayer}
            \begin{pgfonlayer}{back}
                \draw[spedge] (u)--(v);
                \draw[spedge] (x)--(u)--(y);
                \draw[spedge] (x)--(v)--(y);
            \end{pgfonlayer}
        \end{tikzpicture}
        \caption{\label{fig:path-obstruction:after}The resulting structure after performing contraction--uncontraction along the path $xuvy$.}
    \end{subfigure}
    \caption{An obstruction to proving an analogue of the contraction--uncontraction lemma.}
\end{figure}
Consider the graph in \Cref{fig:path-obstruction:before}, which has $x,y\in\doms$, $u\in\onemor x$ and $v\in\onemor y$.
The graph in \Cref{fig:path-obstruction:after} is the graph obtained by performing contraction--uncontraction along the path $xuvy$.
Note that there are $7$ copies of $P_3$ starting at $x$ and not using $y$ in the former graph, whereas there are only $6$ such copies in the latter graph.
Because of this fact, upon performing contraction--uncontraction along $xuvy$, there may be no way to ``locally reroute'' good paths of the form $\subs\doms\cdots\subs\doms\subs\subs$ which use $u$ or $v$ as part of its terminal $\subs\subs$ edge.
If it is not possible to salvage the contract--uncontract lemma, perhaps the notion of benign tumor graphs can be modified to account for this structure, resulting in this structure being accounted for in the maximum likelihood problem.

Currently, we do not see a path around this (and similar) obstacle(s), but we do expect that one exists.

\bibliographystyle{abbrv}
\bibliography{references}

\end{document}